\documentclass[review]{siamart220329}
\usepackage[utf8]{inputenc}

\usepackage{xspace}
\usepackage{amsmath}
\usepackage{bm}
\usepackage{graphicx}
\usepackage{float}
\usepackage{cancel}
\usepackage{enumitem}
\usepackage{lettrine}
\usepackage{siunitx}
\usepackage{textcomp}
\usepackage{gensymb}
\usepackage{pgfplotstable}
\usepackage{booktabs}
\usepackage{wrapfig}
\usepackage{pgfplots}
\usepackage{placeins}
\usepackage{cases}
\usepackage{caption}
\usepackage{subcaption}
\usepackage{cfdpcommands}
\pgfplotsset{compat=1.17}

\newcommand{\eqnref}[1]{\text{Eq.}~(\ref{eq:#1})}

\newcommand{\afigref}[1]{\text{Fig.}~\ref{fig:#1}}

\newcommand{\tblref}[1]{\text{Table}~\ref{tbl:#1}}

\newcommand{\vol}{\Vcl}
\newcommand{\area}{\Acl}
\newcommand{\icellsurf}{\partial \vol_\ibold}
\newcommand{\order}{\text{Q}}
\newcommand{\proj}{\mathbb{P}} 
\newcommand{\Iden}{\mathbb{I}} 

\DeclareMathOperator*{\argmin}{arg\,min}

\title{A Fourth-Order Embedded Boundary Finite Volume Method for the Unsteady Stokes Equations with Complex Geometries
  \thanks{Submitted to the editors May 20, 2022.
  \funding{This work was supported by the Applied Mathematics Program of the U.S. 
DOE Office of Advanced Scientific Computing Research under contract number DE-AC02-05CH11231.
Some computations used the resources of the National Energy Research Scientific Computing Center, 
a DOE Office of Science User Facility supported by the Office of Science of the 
U.S. Department of Energy under Contract No. DE-AC02-05CH11231.}}}

\def \affilCSU {Colorado State University, Fort Collins, CO, USA}
\def \affilLBNL {Lawrence Berkeley National Lab, Berkeley, CA, USA}
\author{Nathaniel Overton-Katz\thanks{\affilCSU (\email{noverton@colostate.edu})}
    \and Xinfeng Gao\thanks{\affilCSU (\email{Xinfeng.Gao@colostate.edu})}
    \and Stephen Guzik\thanks{\affilCSU (\email{Stephen.Guzik@colostate.edu})}
    \and Oscar Antepara\thanks{\affilLBNL (\email{oantepara@lbl.gov})}
    \and Daniel T. Graves\thanks{\affilLBNL (\email{dtgraves@lbl.gov})}
    \and Hans Johansen\thanks{\affilLBNL (\email{hjohansen@lbl.gov})}}

\begin{document}  
\maketitle

\begin{abstract}
A fourth-order finite volume embedded boundary (EB) method is presented for the unsteady Stokes equations.
The algorithm represents complex geometries on a Cartesian grid using EB, employing a technique to mitigate the ``small cut-cell'' problem without mesh modifications, cell merging, or state redistribution.
Spatial discretizations are based on a weighted least-squares 
    technique that has been extended to fourth-order operators and boundary conditions,
    including an approximate projection to enforce the divergence-free constraint.
Solutions are advanced in time using a fourth-order additive implicit-explicit Runge-Kutta method,
    with the viscous and source terms treated implicitly and explicitly, respectively.
Formal accuracy of the method is demonstrated with several grid convergence studies, and results are shown for an application with a complex bio-inspired material.
The developed method achieves fourth-order accuracy and is stable despite the pervasive small cells arising from complex geometries.
\end{abstract}

\begin{keywords}
High-Order Finite Volume, Embedded Boundary, Stokes Equations
\end{keywords}

\begin{MSCcodes}
 	65M08, 76D07
\end{MSCcodes}


\section*{Notation}
\begin{tabbing}
  XXXXX \= \kill 
  $\Dim$ \> spatial dimension, indexed with $d$ \\
  $\xbrm$ \> location in space, e.g., $(x, y, z)$ \\
  $\vol_{\ibold}$ \> volume of a cell $\ibold$ \\
  $\area_{\fbold}$ \> area of a face $\fbold$ \\
  $\kappa$ \> volume fraction relative to the Cartesian cell \\
  $\normvec$ \> surface unit normal vector \\
  $\vec{\Fbrm}$ \> flux tensor, components $\Fbrm_d$\\
  $\ubrm$ \> velocity vector, components $\urm_d$ \\
  $\ibold$ \> grid indices, e.g., $(i, j, k)$ \\
  $\eboldd$ \> unit vector in direction $d$ \\ 
  $\fbold$ \> face indices \\ 
  $\avg{\cdot}$ \> cell-averaged or face-averaged quantity \\
  $\Dbrm$ \> divergence operator \\
  $\Gbrm$ \> gradient operator \\
  $\Lbrm$ \> Laplacian operator \\
  $\proj$ \> projection operator \\
  $O(h)$ \> order of accuracy proportional to a cell length
\end{tabbing}

%
\section{Introduction}
Finite volume methods (FVMs) can achieve high solution accuracy and efficiency on structured Cartesian grids with a high degree of parallelism~\cite{ChomboDesign}.
High-order methods have successfully been created for Cartesian grids \cite{McCorquodale2011} using
    polynomial reconstructions from structured data.
Despite these advantages, representation of complex geometries on structured grids is a significant challenge.
Technologies such as mapped multi-block methods \cite{McCorquodale2015, Guzik2015} can represent moderately complex geometries using structured grids by combining several blocks of curvilinear grids.
However, mapped multi-block grids struggle to represent rough surfaces, and creating quality grids is often time-consuming~\cite{Slotnick2014}.

\begin{figure}[!htbp]
  \centering
  \includegraphics[width=2.1in]{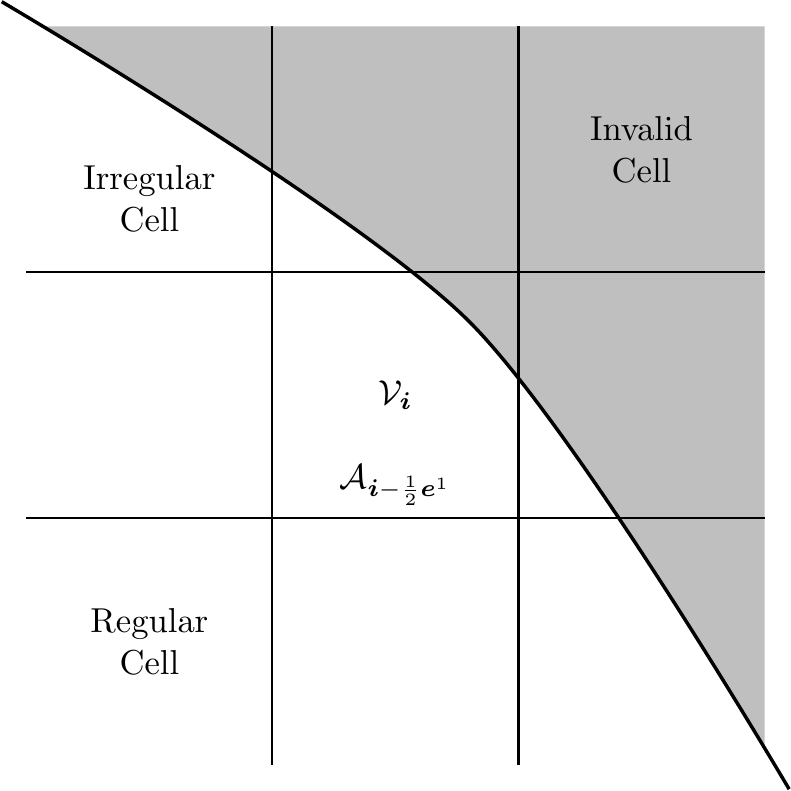}
  \caption{Illustration of a grid for EB methods. The shaded region lies outside the problem domain of interest.}
  \label{fig:ebCells}
\end{figure}

Alternatively, the embedded boundary (EB) method \cite{Aftosmis2000} can represent complicated geometries with little restriction, while maintaining the advantages of structured grid solvers.
EB grids are created by embedding the boundary geometry in a Cartesian grid,
    and cells that intersect the boundaries are cut, as illustrated by \afigref{ebCells}.
The resulting grid is one that is regular away from the boundary, while near the boundary the grid is made up of partial, or ``cut'' cells.
This approach retains the advantages of solving on structured grids on the interior of the domain, while having relatively little restriction in geometries that may be represented.
Additionally, this grid generation process can be done efficiently and quickly in parallel.
However, the presence of cut-cells introduces challenges for both higher-order accuracy and numerical stability.
Specific stencil construction schemes are required near the boundaries, since regular grid-aligned stencils may not be consistent or stable if they use cut-cell values.
In addition, volumes of cut-cells may also become arbitrarily small, and maintaining stability of these small cells requires special care.

While the EB method has been used in a number of complex fluid dynamics applications \cite{Graves2013, Aftosmis2000, Richards2021, Berger2012}, typical applications have only achieved up to second-order accurate solutions, with first-order or 
    inconsistent results near embedded boundaries \cite{Trebotich2015}.
There has been recent progress with improvements in accuracy \cite{Overton2022a}, for example a high-order finite volume EB method for smooth and kinked ($C^0$) domains was demonstrated for Poisson's equation by Devendran et~al.~\cite{Devendran2017}.
The present work extends that fourth-order EB method to solve the time dependent Stokes equations.
Furthermore, a fourth-order additive Runge-Kutta (ARK) scheme is applied to integrate viscous terms implicitly, but any source terms explicitly.
The implicit-explicit (ImEx) time integration is coupled with a projection method \cite{Chorin1968} to enforce the divergence-free constraint on the velocity.

This paper is organized as follows.
\secref{sec:FVMEB} describes the FVM for the EB method, and details the approach taken to solve for the Stokes equations.
Next, in \secref{sec:HOEB}, the algorithm for the fourth-order EB method is developed and applied specifically for the Stokes equations.
The developed algorithm is then verified and validated in \secref{sec:verif} for a number of test cases, and results for complicated geometries are presented and discussed in \secref{sec:results}.
Finally, we draw conclusions and propose directions for future work in \secref{sec:conclusions}

\section{The Finite-Volume Method for Embedded Boundaries}\label{sec:FVMEB}
The FVM discretizes time dependent partial differential equations (PDEs) in ``flux-divergence'' form, 
\begin{equation*}
   \fpartial{}{t}{1} \ubrm(\xbrm, t) + \nabla \cdot \vec{\Fbrm}(\ubrm)= \sbrm(\ubrm,\xbrm,t)\,,
\end{equation*}
where $\ubrm$ is a solution vector that varies in space $\xbrm$ and time $t$, $\vec{\Fbrm}$ is a flux dyad, and $\sbrm$ is a source term
    that may depend on the solution $\ubrm$.
To apply the FVM, the system of PDEs is converted to integral form over volumes $\vol_\ibold$, and the divergence theorem is applied to the flux term $\vec{\Fbrm}$ yielding
\begin{equation}
\label{eq:fvpde}
  \fpartial{}{t}{1}\int_{\vol_\ibold} \ubrm\,\drm\xbrm + \int_{\icellsurf}\vec{\Fbrm}\cdot\normvec\,\drm\xbrm = 
  \int_{\vol_\ibold} \sbrm\,\drm\xbrm\,.
\end{equation}
The volume's surface $\icellsurf$ is split into separate discrete regions $\area_\fbold$, indexed with $\fbold$, so
that $\icellsurf \equiv \cup_\fbold \area_\fbold$.
The surface flux integral then becomes:
\begin{equation*}
  \int_{\icellsurf}\vec{\Fbrm}\cdot\normvec\,\drm\xbrm = \sum_{\area_\fbold \in \icellsurf} \int_{\area_\fbold}\vec{\Fbrm} \cdot \normvec_\fbold\,\drm\xbrm\,,
\end{equation*}
where $\normvec_\fbold$ is the corresponding outward unit normal vector.

FVMs then define averaged quantities on faces and volumes,
\begin{equation*}
  \avg{\ubrm}_\ibold \equiv \frac{1}{\vol_\ibold} \int_{\vol_\ibold} \ubrm\,\drm\xbrm\,,
  \qquad
  \avg{{\Fbrm}}_\fbold \equiv \frac{1}{\area_\fbold} \int_{\area_\fbold}\vec{\Fbrm} \cdot \normvec_\fbold \,\drm\xbrm\,,
\end{equation*}
    and express \eqnref{fvpde} in terms of averages in the discrete ODE form:
\begin{equation}
  \frac{\drm}{\drm t} \avg{\ubrm}_\ibold +  \sum_{\fbold\in\icellsurf} \frac{\area_\fbold}{\vol_\ibold} \avg{{\Fbrm}}_\fbold = \avg{\sbrm}_\ibold\,.
  \label{eq:genFV}
\end{equation}

Note that \eqnref{genFV} has a term of $\vol_\ibold^{-1}$ that must be carefully balanced to avoid numerical stability issues.
As a measure of how small cells are, we define the volume fraction $\kappa$ so that $\kappa \dx^\Dim = \vol_\ibold$\,,
where $\dx$ is the Cartesian grid spacing.
Because cut-cells can be arbitrarily small in our discretization, 
    \eqnref{genFV} must be well-defined in terms of how $\area_\fbold$ and 
    $\avg{{\Fbrm}}_\fbold$ are evaluated on a given small cell as $\vol_\ibold$ (and thus $\kappa$) approach zero.
Regardless, this remains an exact formulation with no approximations, providing the averages and geometric quantities are exact.
In practice, the fluxes and time derivatives require numerical approximations, which characterize the accuracy and stability of the method.
\subsection{Projection Form of the Unsteady Stokes Equations}
The unsteady Stokes equations with constant density are given by
\begin{align}
  \fpartial{}{t}{1} \ubrm
    &= -\nabla p + \nu \Delta \ubrm \,, 
\label{eq:INS1} \\
  \nabla \cdot \ubrm &= 0 \,,
  \label{eq:INS}
\end{align}
where $\ubrm$ is the flow velocity, $p$ the pressure, and $\nu$ is the constant kinematic viscosity.
Boundary conditions for inflows prescribe a velocity $\ubrm = \ubrm_\text{in}$, while outflows are specified by $\nabla\ubrm \cdot \normvec = 0$.
Viscous boundary conditions prescribe a boundary velocity $\ubrm = \ubrm_\text{wall}$. 

%
A Hodge projection operator $\mathbb{P}$ has been used in the finite volume literature \cite{Chorin1968, BCG1989} to enforce the divergence-free velocity field constraint.
We can define this acting on any vector field $\wbrm$ that is not divergence free:
\begin{align}
    \proj(\wbrm) &= \vbrm  \\
    \nabla \cdot \vbrm  &= 0 \\
    \proj(\wbrm) & \equiv \left( \Iden - \nabla \Delta^{-1} \nabla \cdot \right) \wbrm  \,,
\end{align}
    where $\vbrm$ is the divergence free component of $\wbrm$, with appropriate boundary conditions.

\subsection{Finite Volume Projection Formulation}
We choose discretizations for each of the spatial operators in Eqs. (\ref{eq:INS1}-\ref{eq:INS}), and write the resulting discrete equations at cell location $\ibold$ as
\begin{align}
  \label{eq:discreteINS}
  \fpartial{}{t}{1} \ubrm_\ibold
    &= - (\Gbrm p)_\ibold + \nu(\Lbrm \ubrm)_\ibold \,, \\
  (\Dbrm \ubrm)_\ibold &= 0 \,,
\end{align}    
    where $\Dbrm, \Gbrm,$ and $\Lbrm$ are fourth-order finite volume 
    approximations of divergence, gradient, and Laplacian terms, respectively 
(note that from here on, we will drop the subscript $\ibold$ for simplicity).
The goal is to discretize these operators so that
\begin{align}
    \Dbrm \ubrm &= \nabla \cdot \ubrm + O(h^4)\,, \label{eq:divTerm} \\
    \Gbrm \ubrm &= \nabla \ubrm + O(h^4)\,, \label{eq:gradTerm} \\
    \Lbrm \ubrm &= \nabla \cdot \nabla \ubrm + O(h^4)\,, \label{eq:viscTerm}
\end{align}
    in the regular interior of the domain, with some potential loss
    of accuracy near boundaries and in cut-cells.

We use co-located cell-average velocity and pressure, for which a traditional marker-and-cell (``MAC'') 
    staggered-grid discretization of the projection is not available,
    so we instead use an \emph{approximate} projection \cite{Martin2008}. 
This implies that instead of a strictly zero discrete divergence, we allow $\ubrm$ to have a divergence 
    that is at the level of the discretization error.
The equivalent discrete projection is 
\begin{equation}
  \label{eq:discreteProj}
    \Pbrm(\wbrm) = \left( \Iden - \Gbrm \Lbrm^{-1} \Dbrm \right) \wbrm  \, ,
\end{equation}
which requires the inversion of the Laplacian operator over the entire domain.
Using this projection operator, our approximation of the unsteady Stokes equations is
\begin{align}
    \frac{\drm}{\drm t} \ubrm &= \Pbrm \left(\nu \Lbrm \ubrm \right) + O(h^4)     \label{eq:INSproj}
    \\
    \Dbrm \ubrm &= O(h^4) \,,
\end{align}
    and we have eliminated the pressure from the time evolution equations.
The procedure for advancing this system in time is to first do an intermediate update for the viscous terms,
    and then apply the projection \eqnref{discreteProj}, 
    which produces an approximately-divergence free solution at the end of each time step.
This is essentially a higher-order accurate version of the projection operator described by Trebotich et al~\cite{Trebotich2015}.
\subsubsection{Projection Formulation for Open Boundaries}
Boundary conditions for the projection operator can be specified for open domain boundaries,
    so that a given velocity, such as that resulting from the viscous terms, 
    has three separate components: $\wbrm = \Gbrm \psi + \vbrm + \Gbrm \phi$.
First, $\psi$ is the scalar potential flow solution which satisfies only the boundary conditions and
    is thus divergence-free inside the domain.
The two other parts satisfy homogeneous boundary conditions: $\phi$ such that $\Gbrm \phi \cdot \normvec = 0$,
    and the divergence-free part, $\vbrm$, such that $\Dbrm \vbrm \approx 0$ and $\vbrm \cdot \normvec = 0$.
Each of these components is determined from $\wbrm$ by solving the equations:
\begin{align}
\label{eq:projBCs}
    \Lbrm \psi & =  0, \quad \Gbrm \psi \cdot \normvec = \ubrm \cdot \normvec 
    & \hbox{   (potential flow with BCs),} \\
    \Lbrm \phi & =  \Dbrm \wbrm, \quad \Gbrm \phi \cdot \normvec = 0 
    & \hbox{   (interior gradient),} \label{eq:projGradComp} \\
    \ubrm & = \vbrm + \Gbrm \psi = \wbrm - \Gbrm \phi.
\end{align}
This allows $\ubrm$ to be approximately divergence-free and satisfy the correct domain boundary conditions.

Some modifications to this procedure are required for outflow~\cite{Trebotich2015}, given that flow back into
    the domain may occur.
For the divergence operator on the right-hand side of \eqnref{projGradComp} the boundary conditions are:
    $\wbrm$ matches the velocity at inflow $\wbrm = \ubrm_\text{in}$, solid walls require no normal-flow $\wbrm \cdot \normvec = 0$, and outflow boundaries use no boundary condition.
\eqnref{projBCs} and \eqnref{projGradComp} use outflow boundary conditions $\psi = 0$ and $\phi = 0$, respectively.

%
\section{Embedded Boundary Spatial Discretization}
\label{sec:HOEB}
In the embedded boundary approach, there are three categories of cells: regular, irregular, and invalid. 
This distinction between cell types is illustrated in \afigref{ebCells}.
Regular cells are those that are full Cartesian cells, and do not contain a portion of the boundary. 
Irregular or cut-cells, are those which are partial cells because they intersect with the boundary geometry. 
Invalid cells, as the name indicates, are cells that fall outside the domain boundaries and are thus not in the solution domain. 
To denote the EB regions, we intersect the irregular domain, $\Omega$, and $\Upsilon_\ibold$, any regular cell, to denote a particular cell by $\vol_\ibold = \Upsilon_\ibold \cap \Omega$. 

The challenge of EB methods is to approximate the flux terms $\avg{{\Fbrm}}_\fbold$ when regular grid stencils can not be used due to nearby cut-cells.
A general reconstruction evaluates local polynomials and their derivatives on faces to calculate the face-average fluxes in \eqnref{genFV}.
For a structured grid, reconstructed polynomials lead to grid-aligned regular stencils.
When using an EB method, stencils near the boundary depend on the local geometry.
In those cases, we use a weighted least-squares polynomial approximation in a local region of neighboring cells. 
Theoretically, this can produce any order spatial discretization (see Devendran et al.\cite{Devendran2017}), but demonstrating a fourth-order method is the focus of this paper. 
In the following section, the necessary operators needed for the fourth-order EB method are described.

\subsection{Multi-Dimensional Taylor Expansion}
Stencils for the high-order EB method are produced from a multi-dimensional polynomial defined as
\begin{equation*}
  (\xbrm - \bar{\xbrm})^\qbrm = \prod ^{\Dim}_{d=1} (x_d - \bar{x}_d)^{\qrm_d}\,,
  \qquad
  \qbrm! = \prod ^\Dim_{d=1}\qrm_d! \,,
  \qquad
  \left| \qbrm \right| = \sum^\Dim_{d=1}\qrm_d \,,
\end{equation*}
where $\qbrm$ is a \emph{multi-index} or $\Dim$-dimensional non-negative integer vector, and $\bar{\xbrm}$ is a given point in space, which is the center of the interpolation and different for each cell or face. 
For a sufficiently smooth scalar function $\phi$, its multi-dimensional Taylor series of order $\order$ can be written as
\begin{equation}
  \phi(\xbrm) = \sum_{|\qbrm| < \order} \frac{1}{\qbrm!}\phi^{(\qbrm)}(\bar{\xbrm})(\xbrm - \bar{\xbrm})^{\qbrm} + O(h^{\order}) \,,
\end{equation}
where $\phi^{(\qbrm)}$ is the multi-index partial derivative notation,
\begin{equation*}
  \phi^{(\qbrm)}(\xbrm) = \left( \prod ^{\Dim}_{d=1} 
    \frac{\partial^{\qrm_d}}{\partial x_d^{\qrm_d}} \right) \phi(\xbrm) \,,
\end{equation*}
    and any $\qrm_d=0$ implies no derivative.
To fit a multi-dimensional polynomial to cell-averaged data, integration over (regular or irregular) cells is needed.
Using moments to define integration of basis polynomials over a region, we define volume moments and face moments as
\begin{gather}
  m^{\qbrm}_\ibold(\bar{\xbrm}) = \int_{\vol_\ibold} (\xbrm - \bar{\xbrm})^{\qbrm} \,\drm\xbrm \,,
    \qquad
  m^{\qbrm}_\fbold(\bar{\xbrm}) = \int_{\area_\fbold} (\xbrm - \bar{\xbrm})^{\qbrm} \,\drm\xbrm \,.
\end{gather}

\subsection{Flux Reconstruction}
To reconstruct a high-order solution from cell-averaged data, we use a Taylor expansion about cell centers $\bar{\xbrm}_\ibold$ of each uncut-cell $\Upsilon_\ibold$.
For velocity component $\urm_d$, we can define multi-dimensional polynomial coefficients,
$c_{d\ibold}^\qbrm = \frac{1}{\qbrm!}\urm_d^{(\qbrm)}(\bar{\xbrm}_\ibold)$, and approximate each neighboring $\jbold$ cell-average as
\begin{align}
  \avg{\urm_d}_\jbold
  &=
  \frac{1}{\vol_\jbold} \int_{\vol_\jbold} \urm_d \, d \xbrm \nonumber \\
  &=
  \frac{1}{\vol_\jbold} \int_{\vol_\jbold} \sum_{|\qbrm| < \order} \frac{1}{\qbrm!} \urm_d^{(\qbrm)}(\bar{\xbrm}_\ibold)(\xbrm - \bar{\xbrm}_\ibold)^{\qbrm} + O(h^{\order}) \, d \xbrm \nonumber \\
  &=
  \frac{1}{\vol_\jbold} \sum_{|\qbrm| < \order} c_{d\ibold}^\qbrm m^{\qbrm}_\jbold(\bar{\xbrm}_\ibold) + O(h^{\order}) \,.
\label{eq:moments}
\end{align}
The number of coefficients, $c_{d\ibold}^\qbrm$, is $N_p = \frac{(\Dim+\order)!}{\Dim! \,\order!}$.
Determining these $N_p$ coefficients requires solving a linear system of equations with at least that many linearly independent values $\avg{\urm_d}_\jbold$.
We select more cells than coefficients, to establish an \emph{over-determined} least-squares system, which provides
some robustness when there are small cells present in the interpolation neighborhood.

Adapting linear algebra notation, we write the vector $u$ of neighboring cell-average velocities $\avg{\urm_d}_\jbold$,
so our approximation is
\begin{equation*}
  u \approx M \, c \,,
\end{equation*}
where the geometric moment matrix $M$ comes from \eqnref{moments} for each $\jbold$,
    and the vector $c$ contains the multi-dimensional polynomial coefficients $c_{d\ibold}^\qbrm$. 

If we add a diagonal weighting matrix $W$, making a
    \emph{weighted least-squares} (WLS) algorithm, we may improve stability of the stencils 
    (as in \cite{Devendran2017}); further discussion is in section \ref{subsec:weighting}. 
With the weighting matrix, this yields the system to determine coefficients as
\begin{equation}
 c = \argmin_{\tilde{c}} \left\lVert W u - W M \tilde{c} \right\lVert_2  
    \rightarrow c = (W \, M)^{\dagger} \, W \, u \,,
\label{eq:WLS}
\end{equation}
where $(W\,M)^\dagger$ indicates the pseudo-inverse of $(W\,M)$. 

%
\subsubsection{Higher-order EB Viscous Flux Stencil}
\label{sec:ebviscous}
We can use the WLS approach to calculate face-average fluxes for \eqnref{genFV} from cell-average quantities
(as in \cite{Devendran2017}).
This is accomplished by determining the stencil $s_\fbold$ applied to neighbor values $u$,
    derived from the coefficient vector $c$:
\begin{align*}
  \area_\fbold \avg{\Frm_d}_\fbold
  &\equiv s^\intercal_\fbold \, u \\
  &= b^\intercal \, c \\
  &= \left( b^\intercal \, (W \, M)^{\dagger} \, W \right) u \, ,
\end{align*}
    where $b$ is a vector that approximates any face-average flux (\emph{or other quantity}) from the known coefficients.
For a linear flux, such as the viscous term $\Frm_d = \nu \nabla \urm_d$, we determine $b$ (and thus $s_\fbold$) using a Taylor expansion from the 
    face center $\bar{\xbrm}_\fbold$:
\begin{align*}
  \area_\fbold \avg{\Frm_d}_\fbold
  &=
  \int_{\area_\fbold} \nabla \urm_d \cdot \normvec_\fbold \,\drm \xbrm + O(h^\Qrm) \\
  &=
  \int_{\area_\fbold} \nabla \bigg( \sum_{|\qbrm| < \order} c_d^\qbrm (\xbrm - \bar{\xbrm}_\fbold)^{\qbrm} \bigg) \cdot \normvec_\fbold \,\drm \xbrm \\
  &\equiv \sum_{|\qbrm| < \order} c_d^\qbrm b^{\qbrm}_{\fbold} = b^\intercal c \,.
\end{align*}

Each $b^{\qbrm}_{\fbold}$ can be expressed as a sum of \emph{normal-weighted} face moments $m_{\fbold,d}^{\qbrm}$:
\begin{equation}
\label{eq:gradSten}
  b^{\qbrm}_{\fbold} 
  = \sum_{d=1}^{\Dim} \int_{\area_\fbold} \frac{\partial}{\partial{x_d}} 
  (\xbrm - \bar{\xbrm}_\fbold)^{\qbrm} \, \normvec_{\fbold,d} \,\drm \xbrm
  = \sum_{d=1}^{\Dim} \qrm_d m_{\fbold,d}^{\qbrm - \eboldd}
  \,,
\end{equation}
    where $\qbold - \ebold^d$ is required to be positive (derivatives of constants are zero).

For grid-aligned cell faces, normals $\ebold^d$ are constant, so the normal-weighted moments are simply the face moments.
However for curved embedded boundaries, all the components of the normal play a role in the 
    interpolation, especially in the case of inhomogeneous boundary conditions.
Ultimately, the viscous flux stencil $s_\fbold$ depends only on the local neighbor volume and boundary moments, 
    so it may be initialized once per geometry and stored as a sparse matrix operator over the (much smaller) subset of irregular cells.
\subsubsection{Approximate Projection}
The higher-order projection operator starts with cell-average
    velocities, $\avg{\ubrm}_\ibold$, and modifies them to be \emph{approximately}
    divergence-free (\eqnref{INSproj}).
To accomplish this, we require discretizations and boundary conditions for each operator
    in \eqnref{discreteProj}.
First, the Laplacian operator $\Lbrm$ is essentially the same as the
    viscous flux in section \ref{sec:ebviscous}, as a divergence of face-averaged quantities,
    but with different boundary conditions based on \eqnref{projBCs}.
For the cell-average gradient operator, $\Gbrm$, the WLS 
    stencil is similar to \eqnref{gradSten}.
However, in this case we use cell moments instead of face 
    moments to calculate a cell-average gradient:
\begin{equation}
\label{eq:gradStenCell}
  G^{\qbrm}_{d, \ibold} 
  = \int_{\vol_\ibold} \frac{\partial}{\partial{x_d}} 
  (\xbrm - \bar{\xbrm}_\ibold)^{\qbrm} \,\drm \xbrm
  = \qrm_d m_{\ibold}^{\qbrm - \eboldd}
  \,.
\end{equation}    
The gradient operator uses the same boundary conditions as the Laplacian operator.

Finally, we define a cell-average divergence operator $\Dbrm$ using the divergence theorem,
\begin{equation*}
    \int_{\vol_\ibold} \nabla \cdot \ubrm \,\drm \xbrm = \int_{\partial \vol_\ibold} \ubrm \cdot \normvec \,\drm \xbrm \,,
\end{equation*}
so the cell-average of the divergence of the velocity field can be determined from face-average quantities as
\begin{equation*}
    \avg{\nabla \cdot \ubrm}_{\ibold} = \sum_{\fbold\in\icellsurf} \frac{\area_\fbold}{\vol_\ibold} \avg{\ubrm \cdot \normvec}_{\fbold} \,.
\end{equation*}
The divergence flux $\avg{\ubrm \cdot \normvec}_\fbold$ can be constructed similarly to the Laplacian, with 
    an operator that computes a average velocity flux $U_f$ from each velocity component coefficient
\begin{align}
\label{eq:divStenCell}
  U^{\qbrm}_{\fbold} 
  &= \sum_{d=1}^{\Dim} \int_{\area_\fbold} 
  (\xbrm - \bar{\xbrm}_\ibold)^{\qbrm}\normvec_{\fbold,d} \,\drm \xbrm
  = \sum_{d=1}^{\Dim} m_{\fbold,d}^{\qbrm}
  \,.
\end{align}
Again, the terms on regular faces have single component normal vectors $\normvec_\fbold = \ebold^d$, and as a result $\avg{\ubrm \cdot \normvec}_{\fbold} = \avg{\urm_d}_{\fbold}$.
The boundary condition for the divergence is $\ubrm \cdot \normvec = 0$ at any solid
wall, including EB boundaries.
The velocity is specified at inflow boundaries, while outflow has no boundary conditions applied \cite{Colella1999}. 

\subsubsection{Physical Boundary Conditions}
For Dirichlet boundary conditions, each face with a prescribed function uses a boundary average value $\avg{\ubrm}_\fbold= \avg{\ubrm_{\text{bc}}(\xbold)}_\fbold$. A polynomial fitting the function about a face can be reconstructed using:
\begin{align*}
  \area_\fbold \avg{\urm_{\text{bc}, d}}_\fbold
  &=
    \int_{\area_\fbold} \urm_{\text{bc}, d} \,\drm \xbrm \\
  &=
    \int_{\area_\fbold} \sum_{|\qbrm| < \order} c_{\fbold}^\qbrm (\xbrm - \bar{\xbrm}_\fbold)^{\qbrm} \,\drm \xbrm \\
  &=
    \sum_{|\qbrm| < \order} c_{\fbold}^\qbrm m_{\fbold}^{\qbrm} \,.
\end{align*}
Including this additional equation and value into the stencil system will match the boundary condition in a least-squares sense, with a similar method for Neumann boundaries. 
These extra boundary condition equations are used when any neighboring cell included in the reconstruction contains a portion of the boundary, and so can accommodate different parts of the boundary (such as corners, etc.).
This does require that the boundary conditions and derivatives are known at least to order $\order$.

\subsubsection{Interpolation Neighborhoods}
\begin{figure}[!htbp]
  \centering
    \begin{subfigure}[t]{0.48\textwidth}
    \centering
    \includegraphics[width=0.7\linewidth]{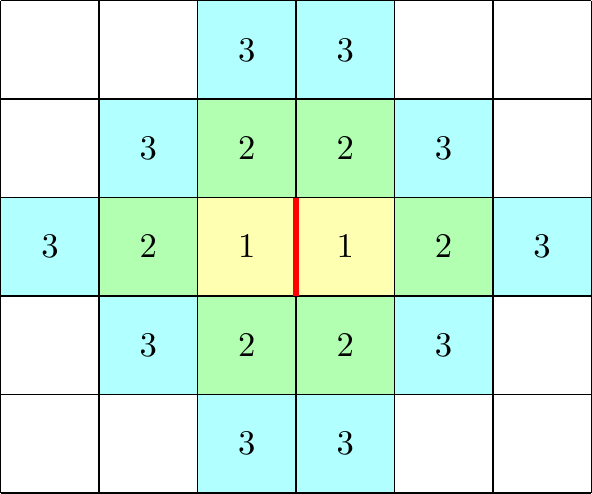}
    \caption{An example regular flux stencil with 18 cells and no boundary conditions.}
  \end{subfigure}
  \,
  \begin{subfigure}[t]{0.48\textwidth}
    \centering
    \includegraphics[width=0.7\linewidth]{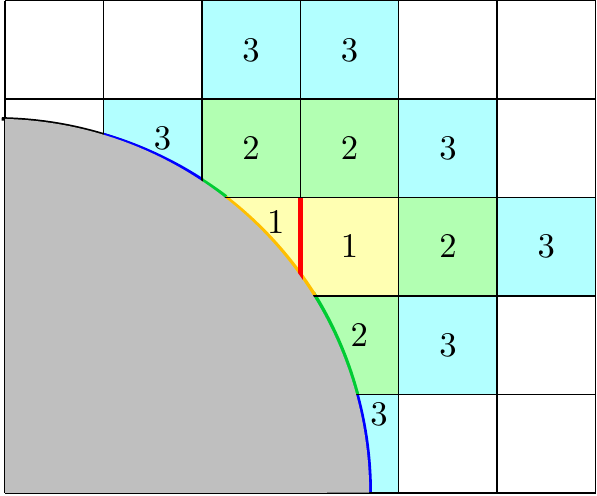}
    \caption{An example irregular flux stencil, using 13 cells and 6 boundary conditions.}
  \end{subfigure}
  \caption{Stencil neighborhoods for flux construction in 2D on the red highlighted faces, with cells numbered according to the Manhattan distance from the indicated red face. Boundary sections included in the stencil are colored in (b).}
  \label{fig:ebStencilsFace}
\end{figure}

\begin{figure}[!htbp]
  \centering
    \begin{subfigure}[t]{0.48\textwidth}
    \centering
    \includegraphics[width=0.7\linewidth]{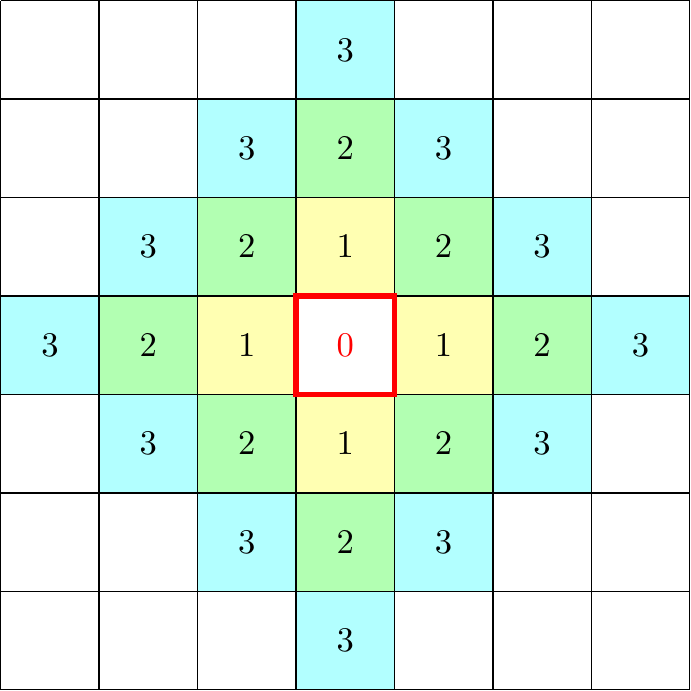}
    \caption{An example regular cell stencil with 25 cells and no boundary conditions.}
  \end{subfigure}
  \,
  \begin{subfigure}[t]{0.48\textwidth}
    \centering
    \includegraphics[width=0.7\linewidth]{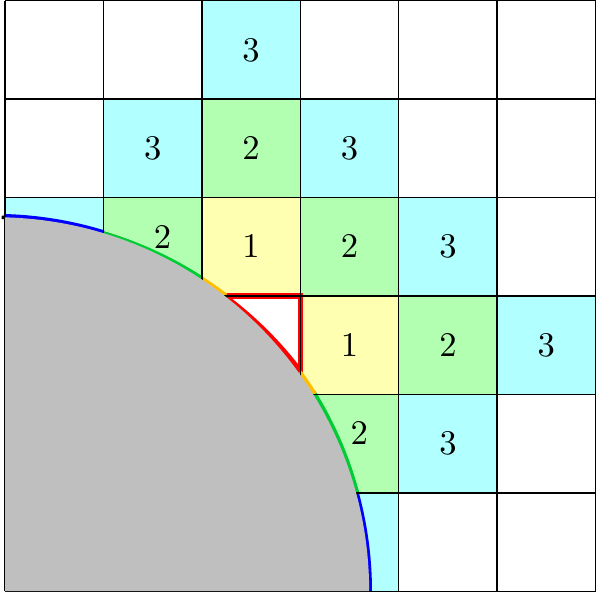}
    \caption{An example irregular cell stencil, using 16 cells and 7 boundary conditions.}
  \end{subfigure}
  \caption{Stencil neighborhoods for cell value construction around the red highlighted cells, with cells numbered according to the Manhattan distance from the red cell. Boundaries included in the stencil are colored in (b).}
  \label{fig:ebStencilsCell}
\end{figure}

An important aspect of the WLS reconstruction is neighborhood selection. Conceptually, we require a large enough
    interpolation neighborhood to attain the desired accuracy and stability properties, while not making it so
    far-reaching so as to create large (expensive), ill-conditioned interpolation matrices.
In this paper, the focus is on fourth-order accurate stencils, which in the general 2D case, requires a 
    minimum of 10 neighbors (with at least $Q$ in each direction) for all of the polynomial coefficients for $|\qbrm| < Q = 4$.
Higher-order derivative operations will require more neighbors to maintain the same accuracy, and values 
    near an embedded boundary will require the boundary condition and further cells because of the lower
    accuracy of one-sided differences.
For example, for a fourth-order WLS viscous operator (Laplacian), flux stencils of 
    radius 3 cells from the reconstructed face can be used for a third-order accurate gradient.
This is illustrated for a radius of 3 in \afigref{ebStencilsFace}(a) for regular regions, and \afigref{ebStencilsFace}(b) in the presence of an EB region with the inclusion of boundary conditions.
However, for stencils that do not specify a boundary condition, such as the divergence on outflow boundaries, this radius of cells may need to be expanded to make the system sufficiently over-determined.
Similarly, for stencils that evaluate cell-averaged values, such as the gradient for the projection operator, cell-average reconstructions are used.
An example of a stencil of radius 3 centered about a cell is shown in \afigref{ebStencilsCell}(a) for regular cells and \afigref{ebStencilsCell}(b) for a cut-cell, following a similar pattern to reconstruction centered on faces.
Cell-average reconstructions may include the cell where the reconstruction takes place, which is indicated as a radius of zero.

\subsection{Weighting and Stability} \label{subsec:weighting}
To reconstruct the stencils required for the spatial discretization in the presence of EB geometries, the WLS method \eqnref{WLS} is used.
First, we select a neighborhood of cells within a given radius (including any cut-cells and their contained subset of the boundary).
Rather than exhaustively searching for an interpolation neighborhood that can determine all the required
    coefficients with finite volume stencils, we use a larger number of neighbors and a weighting scheme that
    assigns relative importance to each cell's entry in the WLS system.
Smaller relative weights mean that cells have less importance and a smaller coefficient in the resulting stencils.
As in Devendran et al. \cite{Devendran2017}, an effective weighting for a fourth-order interpolation is:
\begin{equation}
    W_{\ibold, \fbold} = \max(D_{\ibold, \fbold},1)^{-5} \,,
\end{equation}
where $D_{\ibold,\fbold}$ is the Euclidean distance between cell $\ibold$ center and face $\fbold$ center.
This weighting was shown to improve solution stability and spectral properties, especially when interpolation neighborhoods are large and there are many possible consistent stencils.
\section{Implicit-Explicit Time Marching Method} \label{sec:TimeMarching}
When solving the Stokes equations, the maximum stable time step for the source and viscous terms can be significantly different.
The source terms may generally be non-linear but less stiff, making them well suited for explicit time marching methods.
In contrast, the viscous terms for the Stokes equations are stiff, but linear in each velocity variable, making them well suited for implicit time integration using fast linear solvers.
In the context of EB methods, where cut-cells can become arbitrarily small, the difference of time step stability constraints between terms can be orders of magnitude different.
A hybrid implicit-explicit (ImEx) Runge-Kutta (RK) method \cite{Ascher1997, Zhang2012a} is chosen for the EB algorithm, which allows for the source terms to be updated using an explicit method, and the viscous terms with an implicit method.
This method has successfully been demonstrated for stiff, higher-order finite volume methods for advection-diffusion problems \cite{Zhang2012a}.
To achieve fourth-order accuracy, the six-stage, $L$-stable ImEx scheme ARK4(3)6L[2]SA \cite{Kennedy2003} is used.
At each stage, the implicit time integrator is applied to the viscous term discretization for each velocity component,
    requiring a globally-coupled sparse matrix system to be solved.
To do this efficiently, we use an algebraic multigrid (AMG) solver from PETSC \cite{petsc-user-ref, petsc-efficient}.
After the final stage update from the time integrator, we apply the approximate projection operator
  \eqref{eq:discreteProj} to the predicted velocity to enforce (to the target order of accuracy) 
  the divergence-free constraint \eqref{eq:INS}.
\section{Algorithm Verification}\label{sec:verif}
We use grid convergence to demonstrate the order of accuracy of the algorithm.
The solution error is computed from cell-averages as
\begin{equation*}
  E_\ibold = \avg{\phi(\xbrm, t)}_\ibold - \avg{\phi_{\text{exact}}(\xbrm, t)}_\ibold \,,
\end{equation*}
and convergence rates are evaluated using $L_{\infty}$, $L_2$, and $L_1$ norms computed as
\begin{gather*}
  L_{\infty} = \max_{\ibold \in \Omega} \lvert E_\ibold  \rvert \,, \qquad
  L_{2} = \bigg( \frac{1}{N_c} \sum_{\ibold \in \Omega} E_\ibold^2 \bigg)^{\frac{1}{2}} \,, \qquad
  L_{1} = \frac{1}{N_c} \sum_{\ibold \in \Omega} \lvert E_\ibold \rvert \,,
\end{gather*}
where $N_c$ is the number of cells in the domain $\Omega$.
These solution errors are notably \textit{not} weighted by the cell volumes, in order to properly represent the errors in small cells.
When exact solutions are known, convergence rates are calculated as
\begin{equation}
    \order_{\text{obsv}} = \log_r\left(\frac{\lVert\phi_{h_i} - \phi_{\text{exact}}\rVert}{\lVert\phi_{h_{i+1}} - \phi_{\text{exact}}\rVert}\right) \,.
\end{equation}
where $\phi_{h_i}$ represents the solution at refinement level $h_i$.

The method of manufactured solutions can be used to generate artificial analytic solutions for verification \cite{Salari2000}.
For this method, a manufactured solution of sufficient complexity is chosen, then a source term and boundary conditions are derived which satisfy the governing equations.
However, for the Stokes equations manufactured solutions should still obey the divergence free constraint, limiting the availability of analytical solutions without implementing discrete boundary integral methods.

Algorithm verification of cases without analytic solutions use a standard Richardson extrapolation method to evaluate the convergence by including an additionally refined solution.
Convergence rates with this method \cite{Zhang2012a} are measured by
\begin{equation}
    \order_{\text{obsv}} \approx \log_r\left(\frac{\lVert\phi_{h_i} - \phi_{h_{i+1}}\rVert}{\lVert\phi_{h_{i+1}} - \phi_{h_{i+2}}\rVert}\right)\,.
\end{equation}

We verify the projection and viscous operators separately on simple embedded boundary geometries. 
The projection operator was verified for correctness and stability using the Taylor Green vortex.
The viscous operator with boundary conditions was verified using manufactured solutions in space and time.
The order of accuracy for solving the Stokes equations was verified on a Couette flow between concentric circles.

\subsection{Projection of the Taylor-Green Vortex}
\begin{figure}[!htbp]
  \centering
  \begin{subfigure}[t]{0.48\textwidth}
    \centering
    \includegraphics[width=0.9\linewidth]{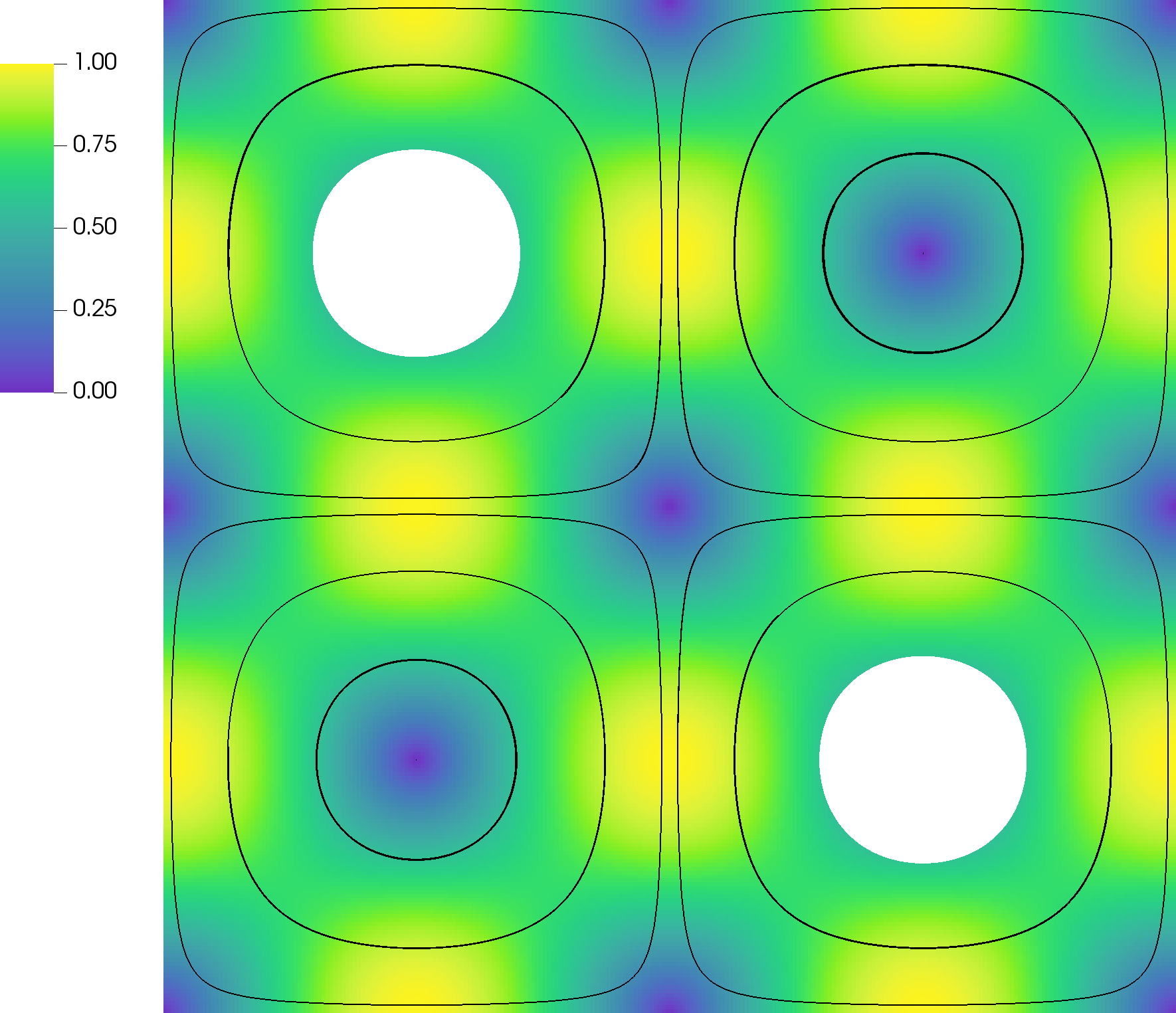}
    \caption{The Taylor-Green velocity magnitude $|\ubrm|$ and streamlines.}
  \end{subfigure}
  \quad
  \begin{subfigure}[t]{0.48\textwidth}
    \centering
    \includegraphics[width=0.9\linewidth]{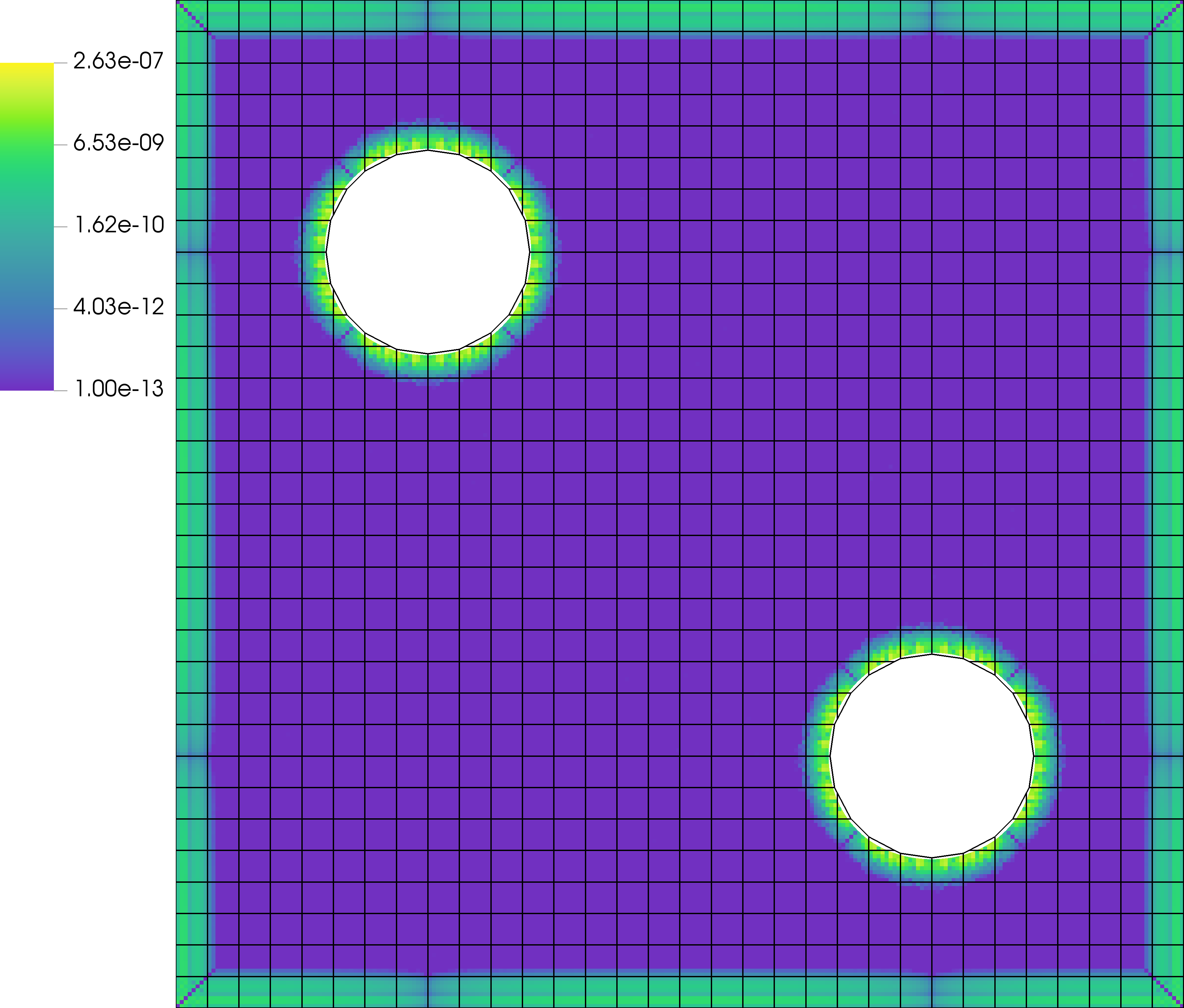}
    \caption{Log plot of the divergence magnitude $\Dbrm \ubrm$ and a coarsened grid representation.}
  \end{subfigure}
  \caption{The Taylor-Green vortex at $h=1/256$.}
  \label{fig:projSolution}
\end{figure}

The approximate projection operator $\Pbrm(\ubrm)$ is demonstrated to be both fourth-order accurate and stable by solving the classic Taylor-Green vortex \cite{Chorin1968} on a unit domain.
The stream function is defined by 
\begin{equation}
\psi = \sin({n \pi x}) \sin({n \pi y}) \,. 
\end{equation}
and the corresponding velocity field in Cartesian coordinates is,
\begin{gather}
    u = \sin({n \pi x}) \cos({n \pi y}) \,, \quad  
    v = \cos({n \pi x}) \sin({n \pi y}) \,, 
\end{gather}
which is analytically divergence free.
For this case, we choose the period $n=2$, and as a result there is no flow through the domain boundary.
In addition, EB boundaries are cut along the contours of the stream function at $\psi = -0.8$.
Initial conditions use the exact velocity profiles, integrated to fourth-order cell-averages, while projection boundary conditions are specified (no-flow, but not viscous, walls).
The geometry and velocity field are shown in \afigref{projSolution}(a).

\begin{figure}[!htbp]
  \centering
  \begin{subfigure}[t]{0.48\textwidth}
    \centering
    \includegraphics[width=0.7\linewidth]{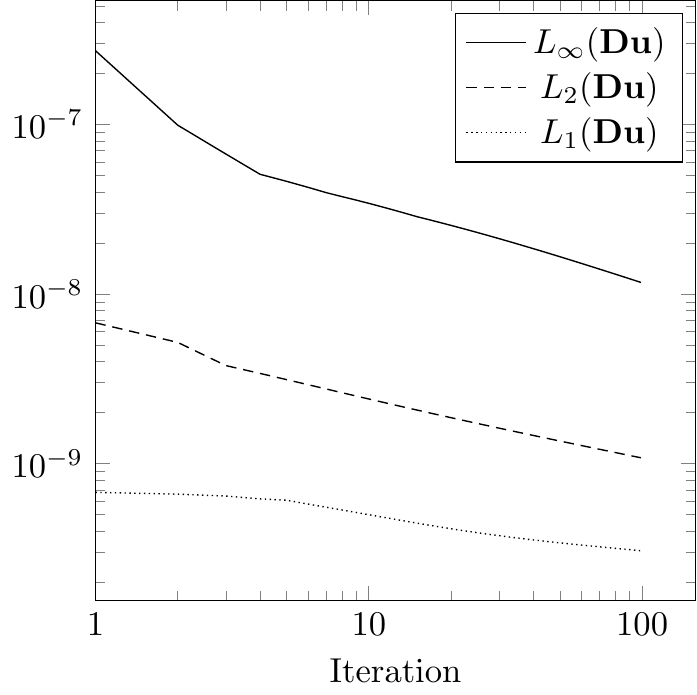}
    \caption{Error norms of the velocity divergence.}
  \end{subfigure}
  \quad
  \begin{subfigure}[t]{0.48\textwidth}
    \centering
    \includegraphics[width=0.7\linewidth]{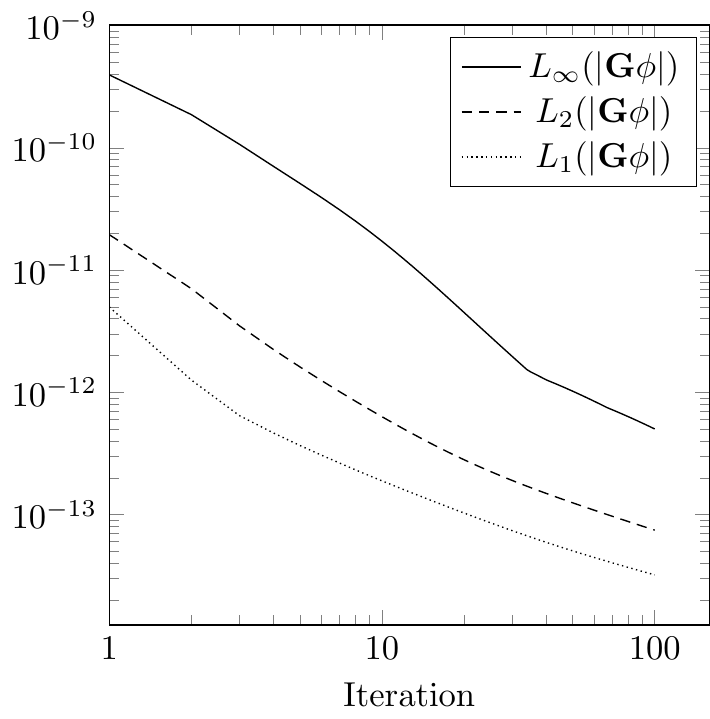}
    \caption{Error norms of the pressure gradient.}
  \end{subfigure}
  \caption{Solution error norms from repeated projections of the Taylor-Green vortex on grid size $h=1/256$.}
  \label{fig:projStability}
\end{figure}

Stability of the approximate projection operator is demonstrated by showing that repeated applications of the projection on a velocity field reduce the discrete divergence monotonically towards zero.
The velocity divergence $\Dbrm \ubrm$ and pressure gradient $\Gbrm \phi$ are computed on a grid with cell size $h=1/256$ for 100 projection applications and plotted in \afigref{projStability}.
These quantities strictly decrease, showing stability of the fourth-order approximate projection.

\begin{figure}[!htbp]
  \centering
  \includegraphics[width=2.3in]{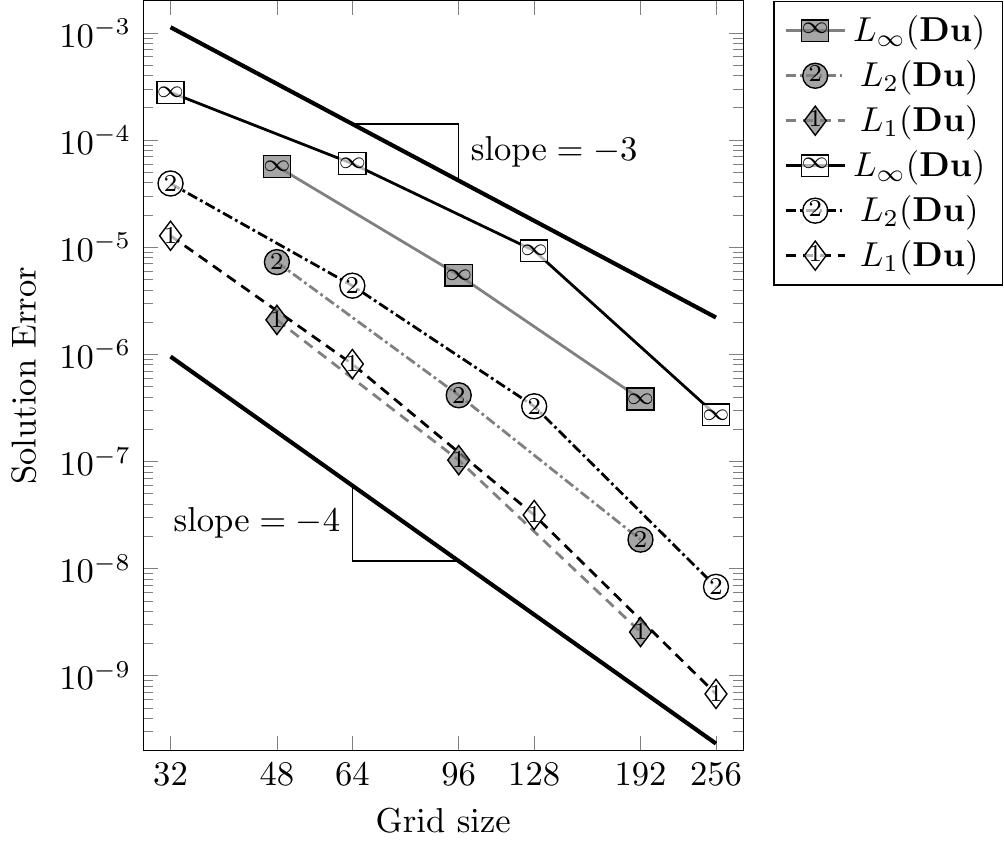}
  \caption{Convergence rates of the divergence for the Taylor-Green vortex in 2D. The gray series of points are coarsened from a fine level of $h=1/192$, and the white series coarsened from $h=1/256$.} 
  \label{fig:projectionConvg2D}
\end{figure}

Convergence tests are performed to ensure the targeted fourth-order accuracy is achieved.
The projection is applied once for grids of decreasing refinement, and the divergence field $\Dbrm \ubrm$ is used to evaluate the solution error.
The finest levels are chosen with cell sizes $h=1/256$ and $h=1/192$, and subsequent coarser levels each double $h$.
The $L_1$, $L_2$, and $L_\infty$ errors are calculated and plotted in \afigref{projectionConvg2D}, where convergence rates reach and exceed expected fourth-order accuracy.
Divergence errors are the largest in cut-cells at the embedded boundaries, as seen in \afigref{projSolution}(b), and dominate the $L_\infty$ errors.
Nevertheless, the solution in cut-cells still converges at the expected rate.

%
\subsection{Manufactured Solution for Diffusion Inside a Circle}
\label{sec:resultsDiff}
To demonstrate the accuracy of the viscous operator in the Stokes equations, we use a manufactured solution and solve without the projection operator.
This completely decouples the velocity equations.
The algorithm targets fourth-order in space and time using the high-order stencils described in section \ref{sec:TimeMarching} and the ImEx scheme in section \ref{sec:HOEB}.
A circular domain is created of radius $0.3$ centered about the point $(0.5, 0.5)$.
The manufactured solution is the same for each component, defined by
\begin{equation} 
    \urm_d(\xbrm, t) = \sin(2\pi t)\sin(R^2 - (\xbrm - \xbrm_0)^2)
    \label{eq:manDiff}
\end{equation}
where $R=0.3$ matches the domain radius, and $\xbrm_0=(0.5, 0.5)$ gives the domain center.
The embedded boundaries are specified by Dirichlet conditions with values determined by the manufactured solution.
Following the method of manufactured solutions, a source term is added to balance the equation where
\begin{align*} 
    \srm_d(\xbrm, t) =& 2\pi\cos(2\pi t)\sin(R^2 - (\xbrm - \xbrm_0)^2) \\
    &\quad -  \sin(2\pi t) \sum_d^\Dim \left(-4 \xrm_d \sin(R^2 - (\xbrm - \xbrm_0)^2) + 2\cos(R^2 - (\xbrm - \xbrm_0)^2) \right)\,.
    \label{eq:mmsSource}
\end{align*}
The source term is evaluated explicitly in time, while the Laplacian term is evaluated implicitly.
The solution is initialized using fourth-order cell-averages at time $0.125$, and a viscosity of $\nu = 1$.
On the finest level, with grid size $h=1/128$, a time step of $\Delta t = 0.1$ is taken to advance the solution forward for 128 steps.
Subsequent coarser levels double the cell size and time step, while halving the number of time steps to reach the same end time.
Using the chosen exact solution in \eqnref{manDiff}, the $L_1$, $L_2$, and $L_\infty$ errors and convergence rates are calculated and compiled in \tblref{diffConverg2D}.
Third-order truncation error is anticipated at the embedded boundaries, and will dominate the $L_\infty$ norm.
However, the $L_1$ and $L_2$ norms still attain fourth-order accuracy because the embedded boundary is only codimension one \cite{Johansen1998}.
Fourth-order accuracy is demonstrated in these error norms, verifying the algorithm.

\begin{table}[!htbp]
  \centering
  \caption{Viscous operator convergence errors and rates for diffusion inside a circle.}
    \begin{tabular}{c|cccccc}
      $N$ & $L_1$ & Rate($L_1$) & $L_2$ & Rate($L_2$) & $L_\infty$ & Rate($L_\infty$) \\ \hline
      16  & 3.676e-07 & & 5.323e-07 & & 1.271e-06 & \\
      32  & 1.421e-08 & 4.693 & 2.111e-08 & 4.656 & 6.810e-08 & 4.223 \\
      64  & 6.449e-10 & 4.463 & 9.529e-10 & 4.470 & 3.186e-09 & 4.418 \\
      128 & 3.688e-11 & 4.128 & 5.467e-11 & 4.123 & 2.195e-10 & 3.860
    \end{tabular}%
  \label{tbl:diffConverg2D}
\end{table}


%
\subsection{Spherical Couette Flow}
\begin{figure}[!htbp]
  \centering
  \begin{subfigure}[t]{0.48\textwidth}
    \centering
    \includegraphics[width=0.95\linewidth]{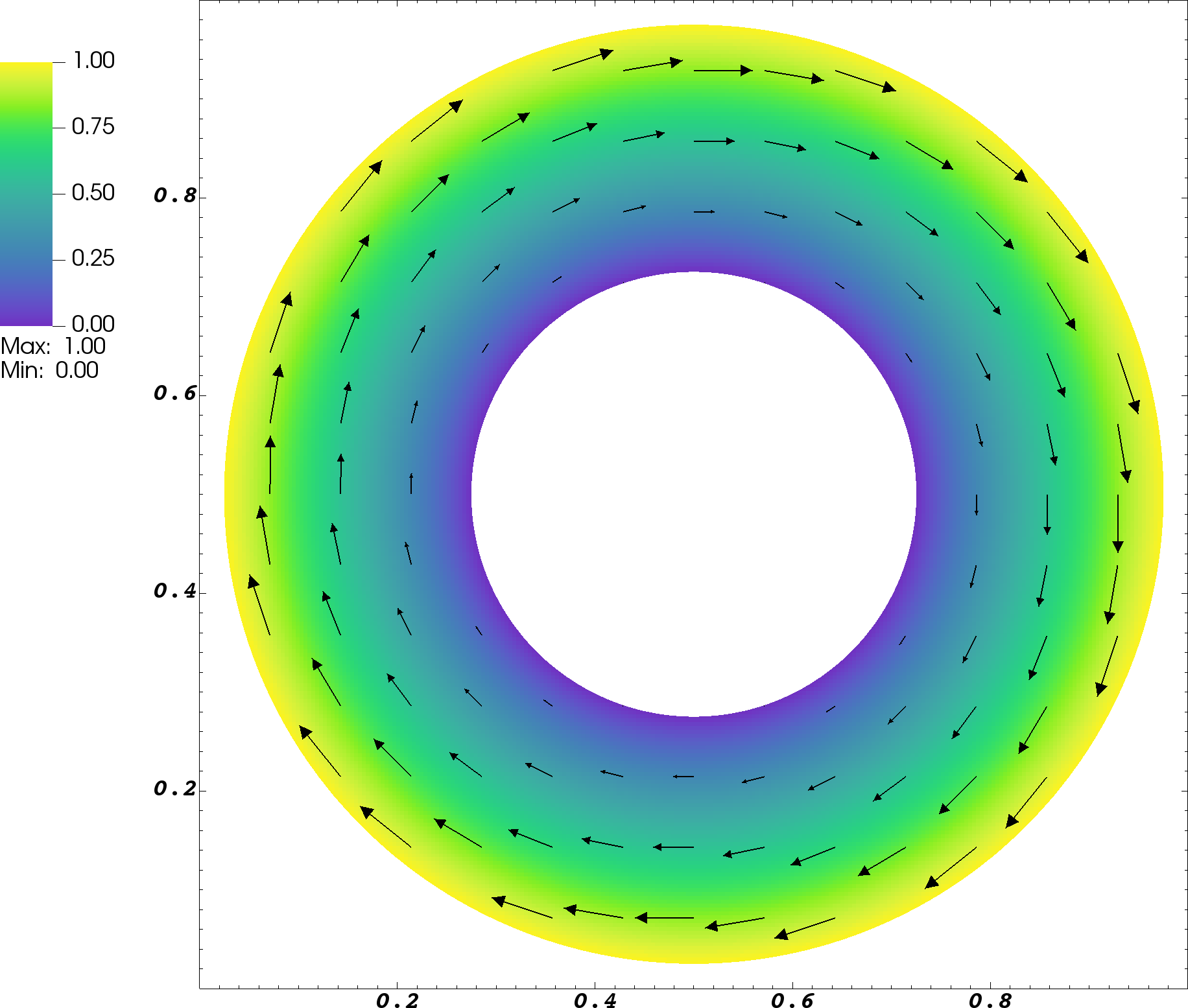}
    \caption{Velocity magnitude and vectors for the $h=1/256$ solution.}
  \end{subfigure}
  \quad
  \begin{subfigure}[t]{0.48\textwidth}
    \centering
    \includegraphics[width=0.85\linewidth]{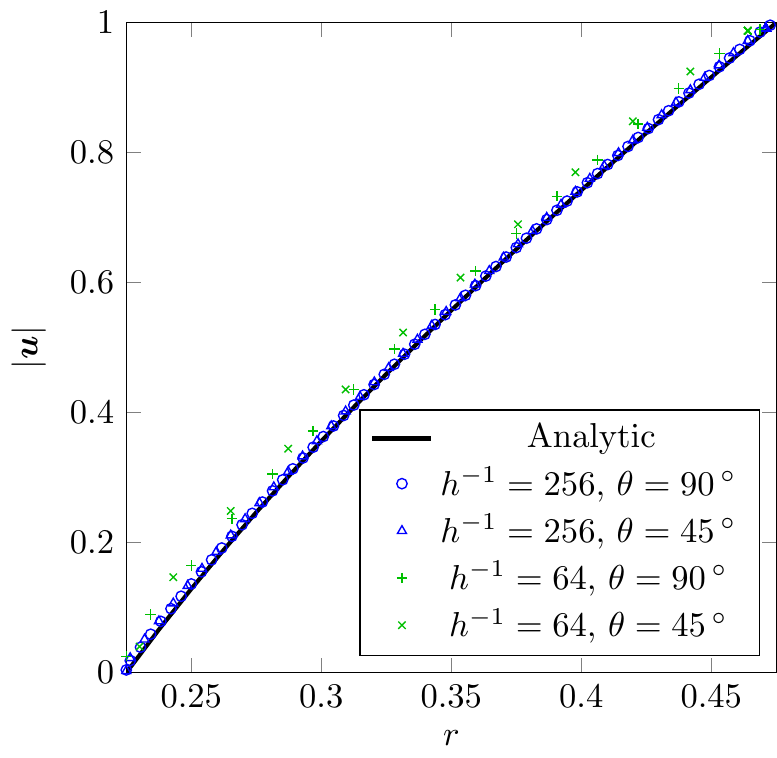}
    \caption{Velocity magnitude profiles compared to the analytic solution.}
  \end{subfigure}
  \caption{The two-dimensional circular Couette flow at steady state.}
  \label{fig:couetteSol}
\end{figure}

\begin{figure}[!htbp]
  \centering
  \includegraphics[width=2.4in]{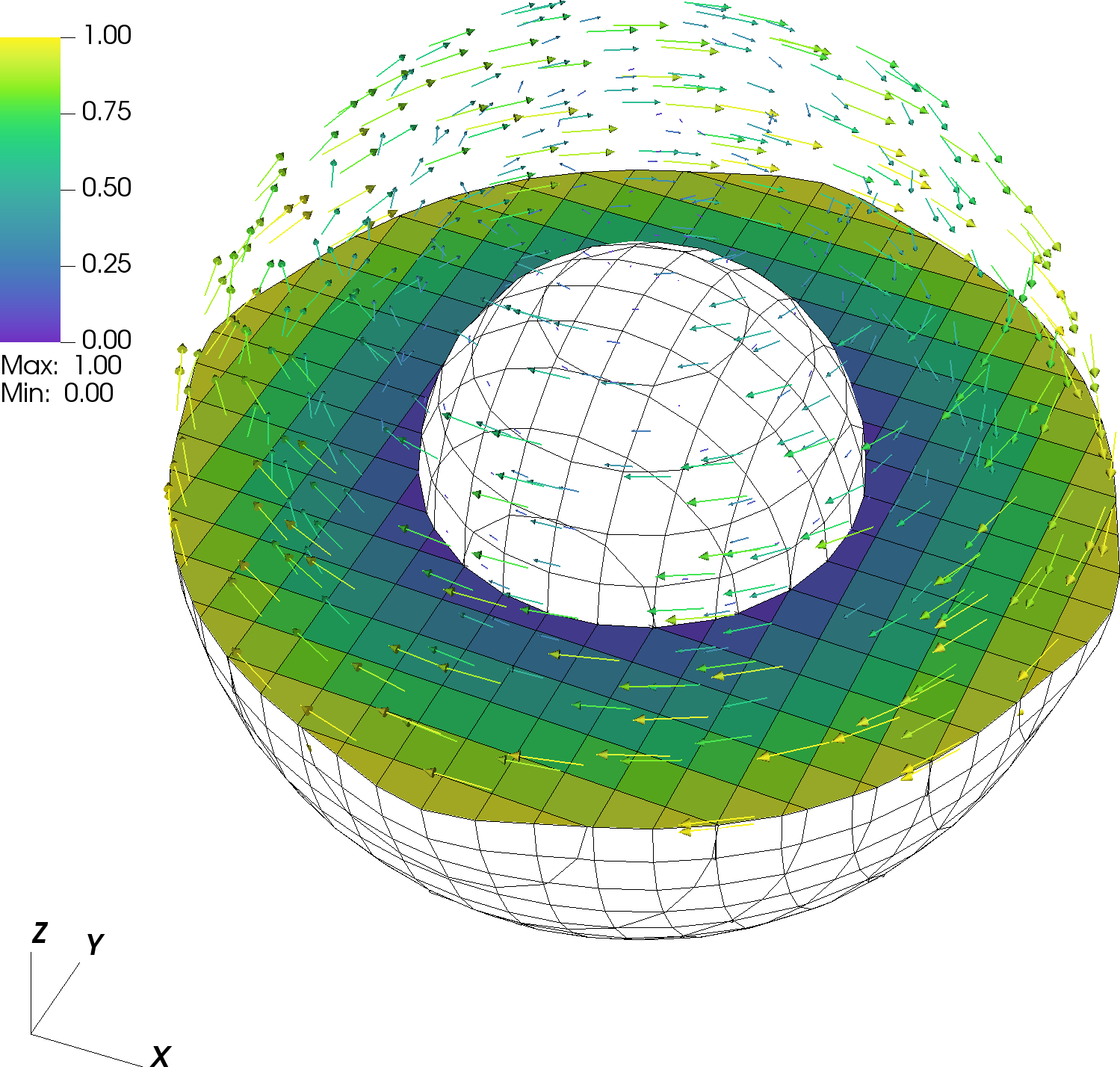}
  \caption{The three-dimensional spherical Couette flow at steady state. The velocity vector field is shown, and a slice along the $x$-$y$ plane shows the grid and velocity magnitude.}
  \label{fig:couetteSol3D}
\end{figure}

A pair of concentric spheres are created with radii $r_{\text{inner}}=0.25$ and $r_{\text{outer}}=0.475$.
The inner sphere is held stationary, while the outer sphere is rotated about the z-axis with constant angular velocity $\omega_\text{outer}=\frac{1}{0.475}$ so that the outer sphere has a peak tangential velocity of 1.
Velocity of the rotating sphere is specified purely tangential to the surface by
\begin{equation*}
    u_\text{tan} = \sin(\varphi)R\omega
\end{equation*}
where $R$ is the sphere radius, $\omega$ the angular velocity, and $\varphi$ the polar angle from the z-axis.
Converting this to Cartesian coordinates prescribes a velocity field
\begin{equation*}
    u = \omega y \,, \qquad
    v = - \omega x \,, \qquad
    w = 0 \,,
\end{equation*}
where $x$ and $y$ are the Cartesian coordinates on the sphere centered at the origin.
The solution is initialized to zero velocity uniformly with a viscosity of $\nu = 1$.
Using the developed method in this work, the Stokes equations are solved to steady state. 
A two-dimensional case normal to the z-axis is shown in \afigref{couetteSol}(a), and the three-dimensional version in \afigref{couetteSol3D}.

\begin{figure}[!htbp]
  \centering
  \begin{subfigure}[t]{0.48\textwidth}
    \centering
    \includegraphics[width=0.9\linewidth]{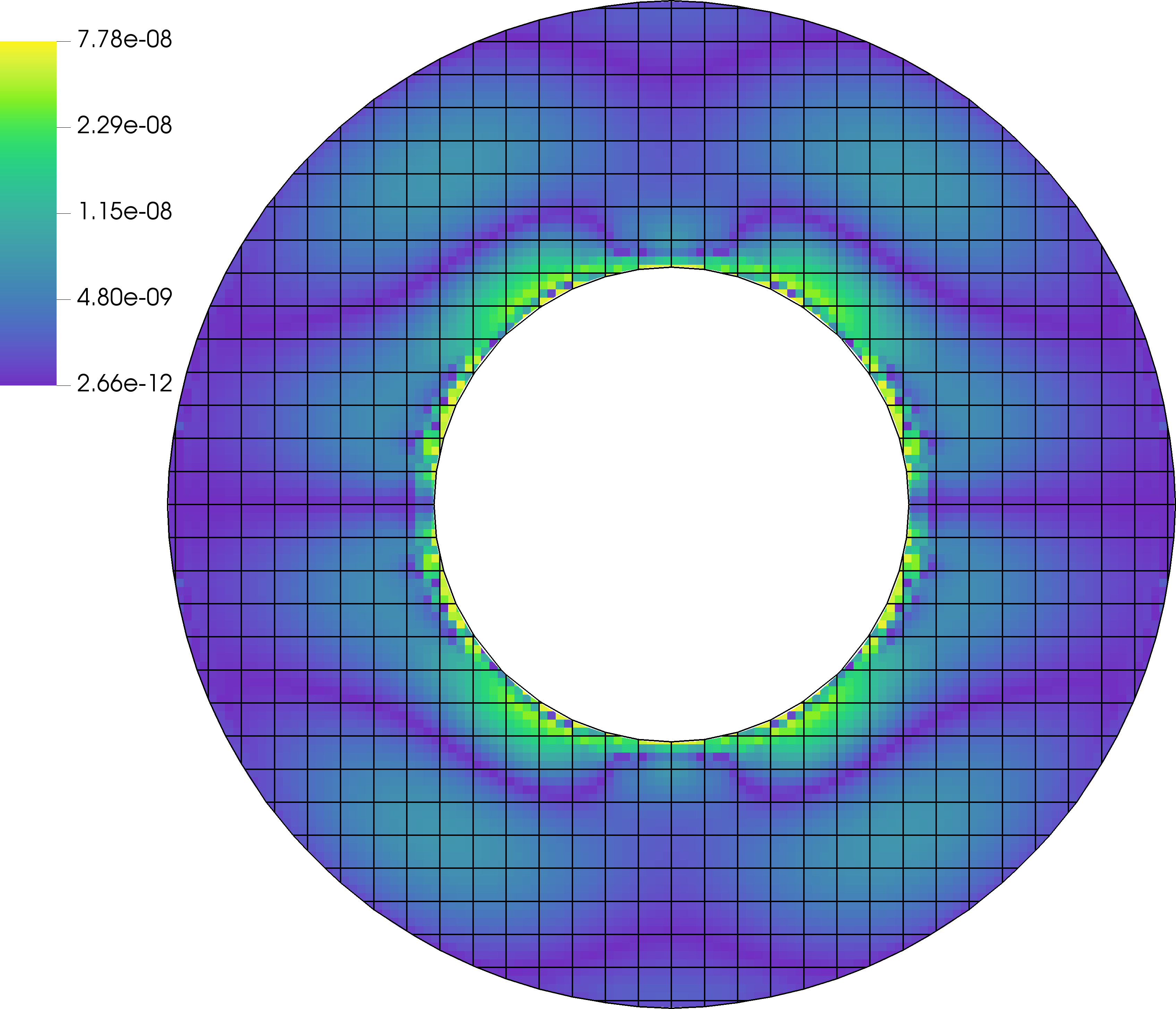}
    \caption{Solution errors of the $u$ velocity on the $h=1/128$ grid, with a coarsened grid representation.}
  \end{subfigure}
  \quad
  \begin{subfigure}[t]{0.48\textwidth}
    \centering
    \includegraphics[width=0.9\linewidth]{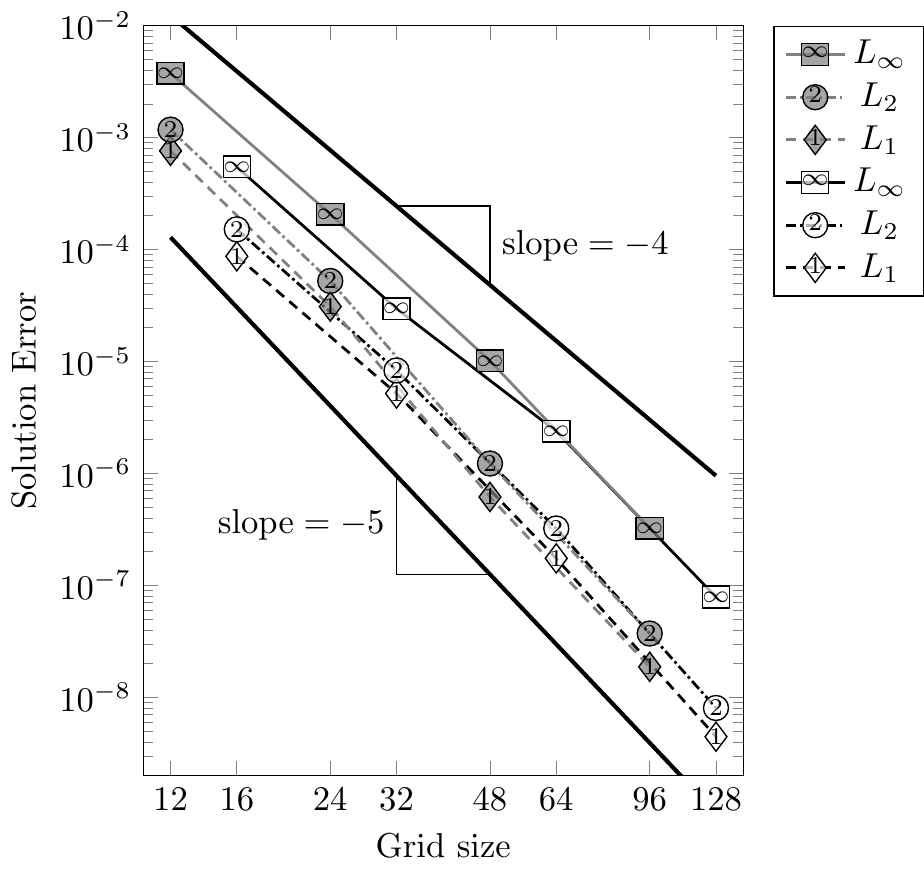}
    \caption{Convergence rates of the $u$ velocity.
    The gray series of points are coarsened from a fine level of $h=1/192$, and the white series from $h=1/256$.}
  \end{subfigure}
  \caption{The two-dimensional circular Couette flow solution errors.}
  \label{fig:couetteErr2D}
\end{figure}

\begin{figure}[!htbp]
  \centering
  \includegraphics[width=0.4\linewidth]{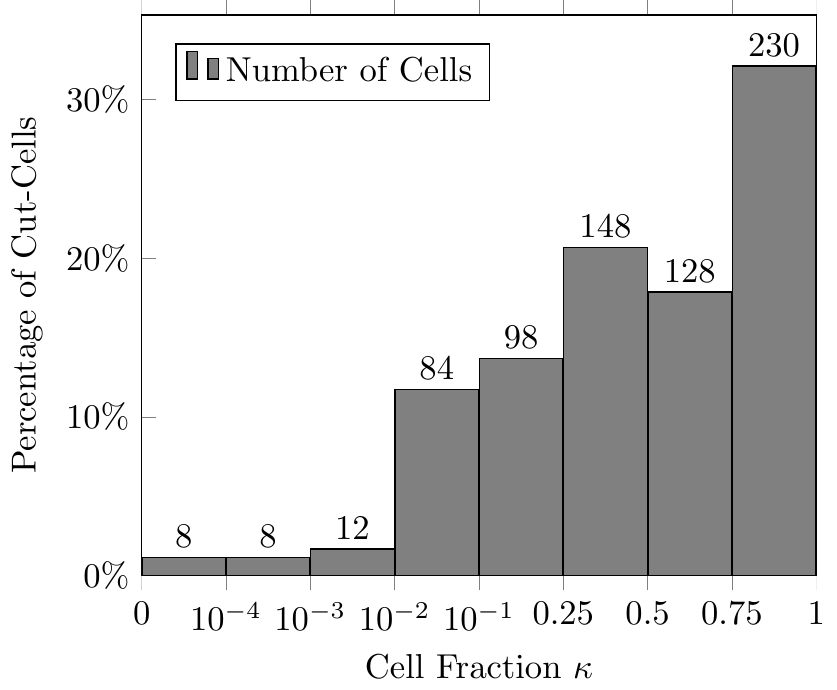}
  \caption{Distribution of cut-cell sizes for the two-dimensional Couette flow at grid size $h=1/256$.} 
  \label{fig:couetteVofFracs}
\end{figure}

On the finest level, of grid size $h=1/256$, a time step of $\Delta t=$ \num{1e-3} is taken to advance the solution forward until the divergence norm converges in the two-dimensional case.
To validate the solutions, radial profiles of the solution in two dimensions are compared to the analytic solution \cite{Masatsuka2013}
\begin{equation*}
    u_\theta(r) = \omega_\text{outer} r_\text{outer} \frac{{r}/{r_\text{inner}} - {r_\text{inner}}/{r}}{{r_\text{outer}}/{r_\text{inner}} - {r_\text{inner}}/{r_\text{outer}}} \,,
\end{equation*}
in \afigref{couetteSol}(b), where good agreement is observed.
Convergence rates are measured and shown in \afigref{couetteErr2D}(b) using the Richardson extrapolation method for a series of grids coarsened by a factor of 2 from refinements of both $h=1/256$ and $h=1/192$.
Solution norms for the $u$ and $v$ velocity components only have differences on the order of machine precision, since the flow is symmetric, and so only one set of errors are shown.
Results show that fourth-order accuracy is achieved or exceeded for all solution norms. 
The $L_\infty$ norm in particular is dominated by cut-cell values, as shown in \afigref{couetteErr2D}(a).
On the finest grid there are \num{39364} valid cells, of which \num{716} are cut-cells with volume fractions as small as $\kappa$ = \num{1.317e-5}, demonstrating the robustness of the method.
The distribution of cut-cell sizes is shown in \afigref{couetteVofFracs}.


%
\section{Results}\label{sec:results}
Using the verified and validated fourth-order EB algorithm, Stokes flow over a circle and sphere in a channel are tested, and convergence order is verified.
The Stokes equations are also used to solve for flow through a complex bio-inspired material of engineering interest.

\subsection{Steady Stokes Flow Over a Sphere in a Channel}
\begin{figure}[!htbp]
    \centering
    \includegraphics[height=1.5in]{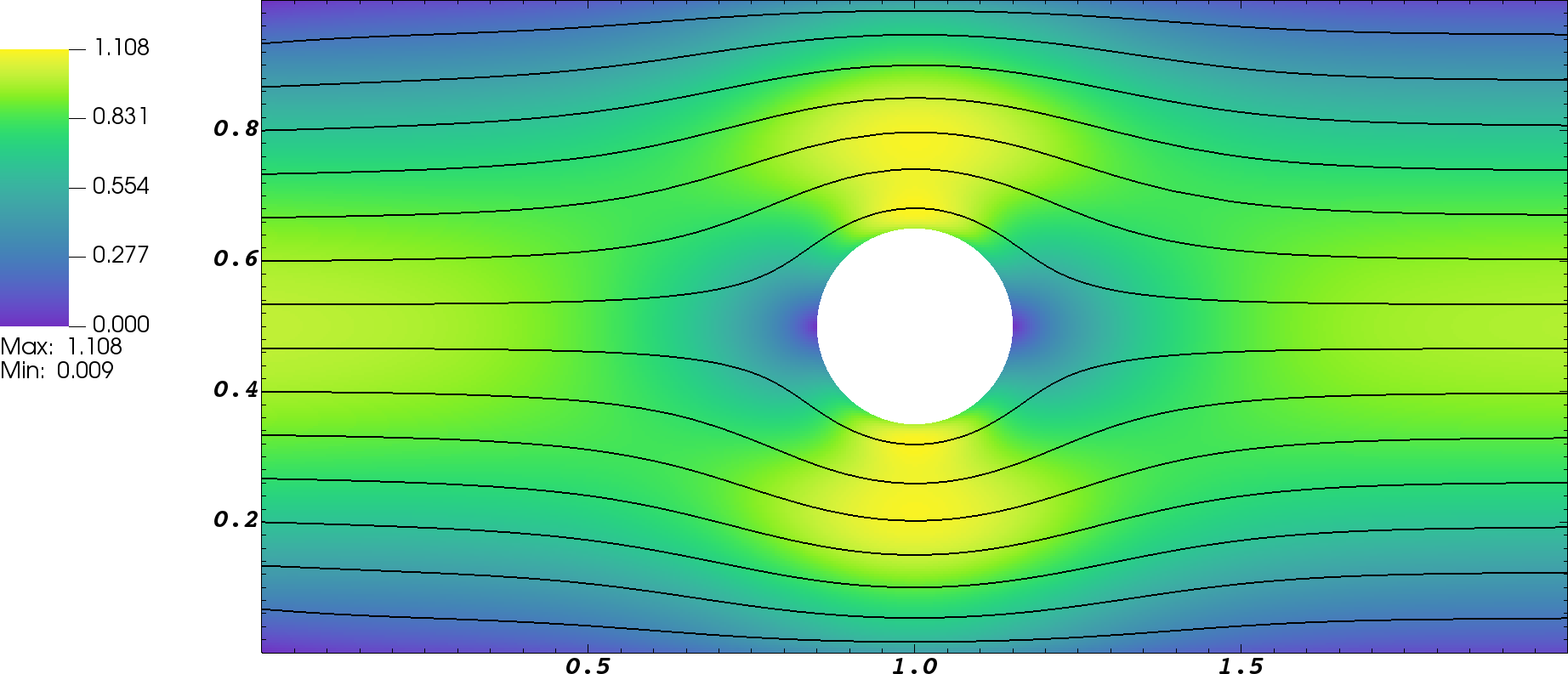}
    \caption{Stokes flow for a circle in a channel, where the left boundary is an inlet, the right an outlet, and all other boundaries walls. Streamlines are shown in black, and the contours plot the velocity magnitude.} 
  \label{fig:channel2D}
\end{figure}

\begin{figure}[!htbp]
    \centering
    \includegraphics[height=2.2in]{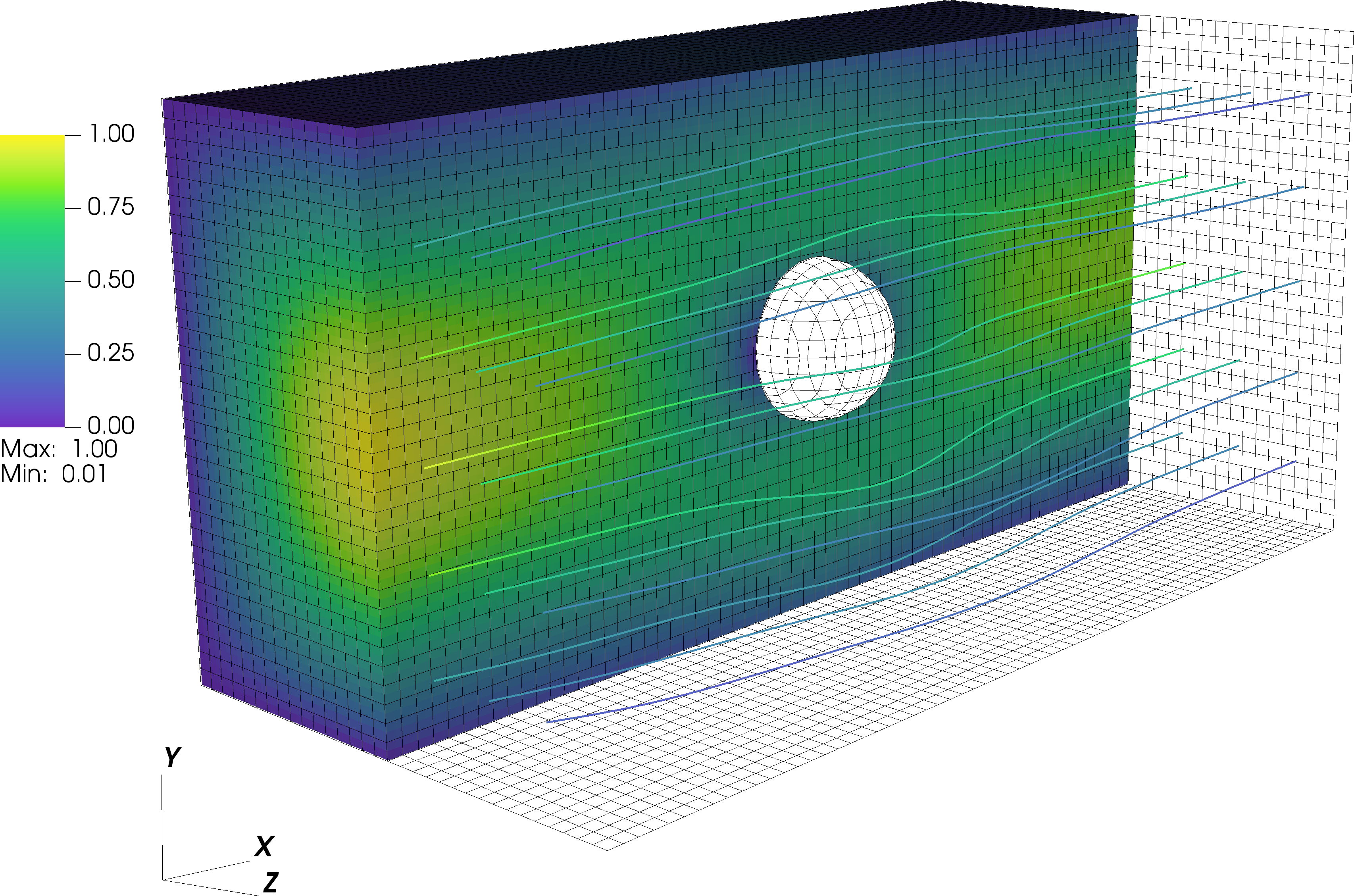}
    \caption{Stokes flow over a sphere in a channel, where the left boundary is an inlet, the right an outlet, and all other boundaries walls. The streamlines and contour plot show the velocity magnitude.}
  \label{fig:channel3D}
\end{figure}

A square channel is generated with a channel length of $2$ in the $x$-direction, and an inlet length of $1$ in the $y$ and $z$ directions. 
A sphere of radius $r=0.15$ is centered in the channel at $x=1$.
A developed inflow of peak value 1 is specified at the inflow on the left most boundary, while an outflow is specified for the right most boundary.
All other boundary conditions are specified as walls.
The Stokes equations are solved to steady state in both two-dimensions, in \afigref{channel2D}, and three-dimensions, in \afigref{channel3D}.
Although no analytic solution for this flow can be used for comparison, qualitatively the flow is observed to be highly symmetric and respect the imposed boundary conditions.
Notably, the inflow boundary and outflow nearly match, and have streamlines normal to the boundaries as expected.
At wall boundaries, solution velocities approach zero, although particularly along the sphere there is some discrepancy due to the projection and viscous operators smearing the third-order boundary solution with the fourth-order interior.

\begin{figure}[!htbp]
  \centering
  \includegraphics[width=2.4in]{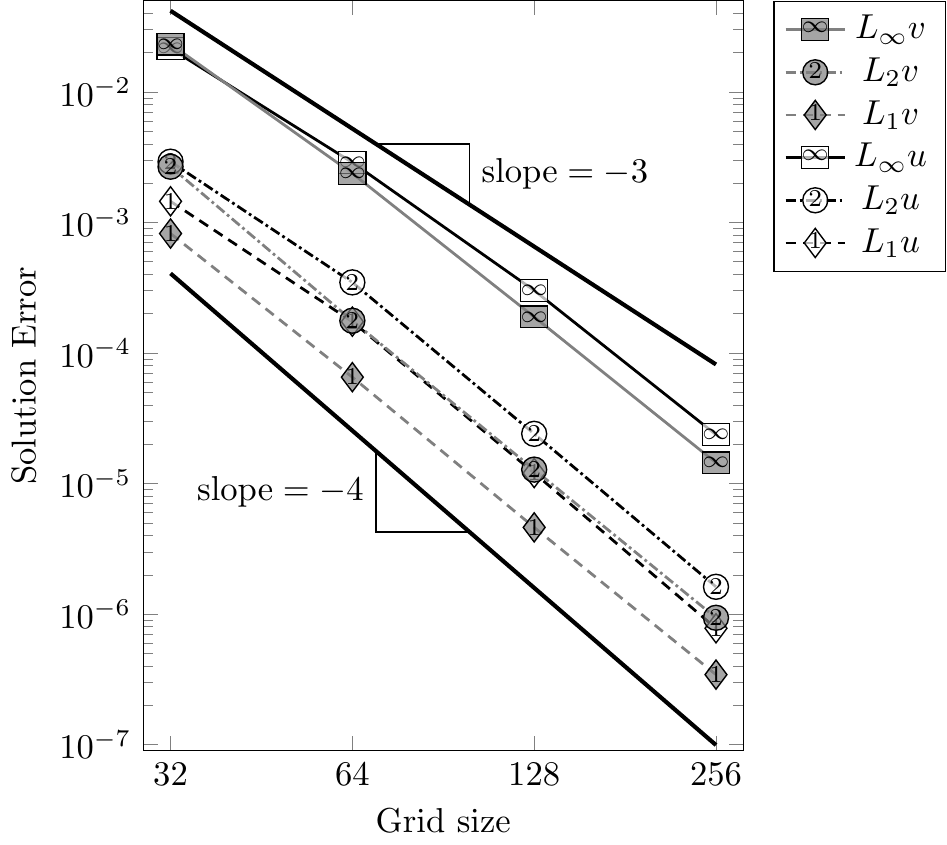}
  \caption{Convergence rates of the 2D circle in a channel.} 
  \label{fig:channelConvg2D}
\end{figure}

\begin{figure}[!htbp]
    \centering
    \includegraphics[height=1.5in]{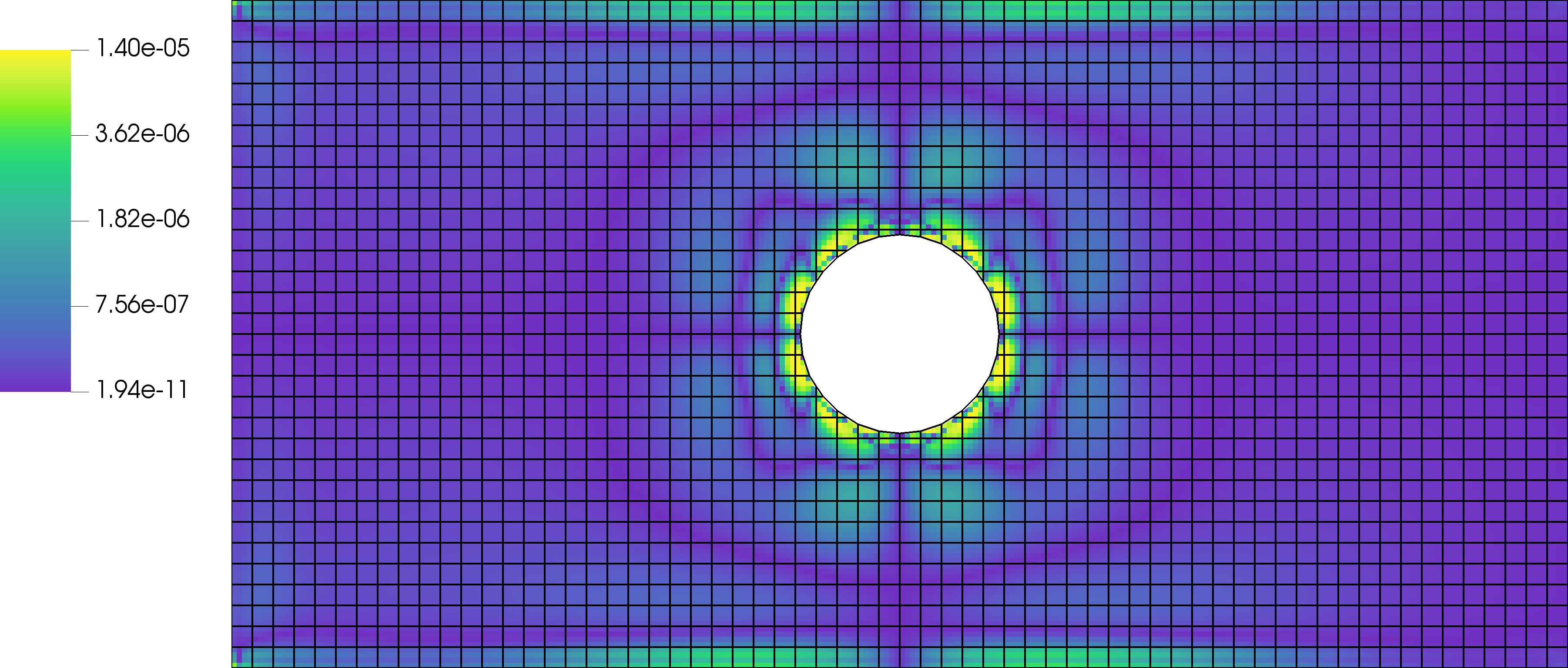}
    \caption{Solution errors for the $u$ velocity of the 2D Stokes flow over a circle in a channel on the $h=1/265$ grid, with a coarsened grid representation.} 
  \label{fig:channel2Derr}
\end{figure}

Convergence rates for the $u$ and $v$ velocity from Richardson extrapolation are plotted in \afigref{channelConvg2D}.
$L_1$ and $L_2$ norms indicate fourth-order accuracy, but the $L_\infty$ norm lags slightly.
The plot of solution error (\afigref{channel2Derr}) indicates lower accuracy in cut-cells and parts of the domain boundary.

\subsection{Stokes Flow for Bone Scaffolding}
\begin{figure}[!htbp]
    \includegraphics[width=\linewidth]{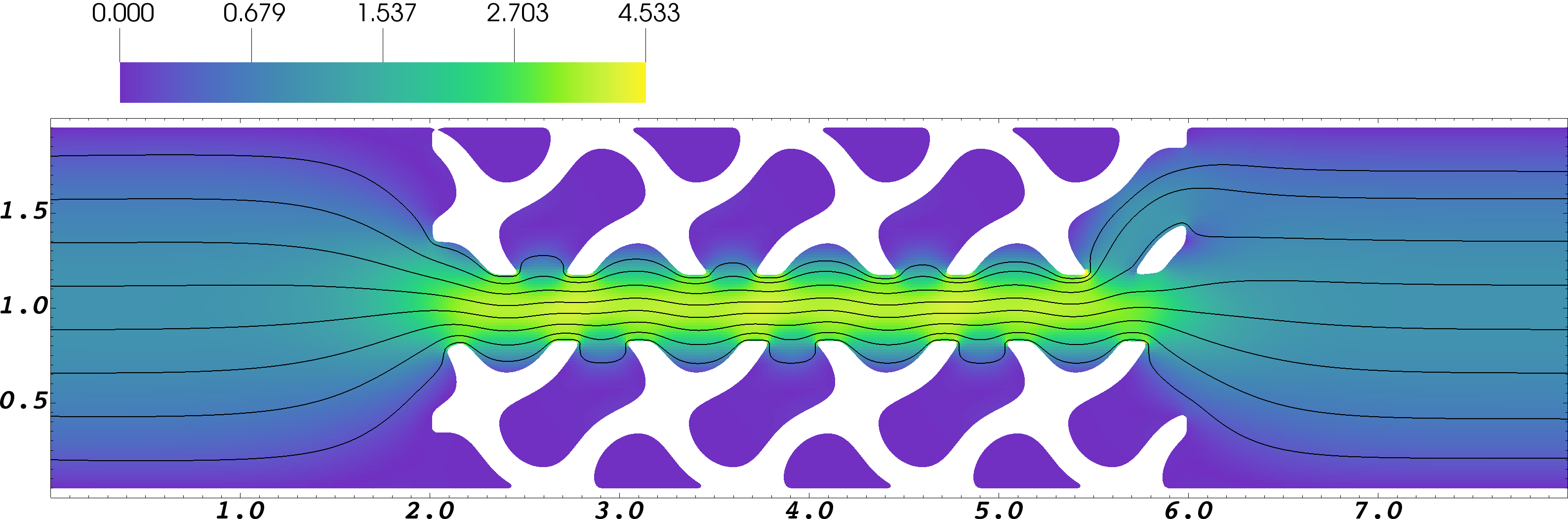}
    \caption{Velocity magnitude and streamlines for a two-dimensional representation of the bone scaffold geometry.}
    \label{fig:bone}
\end{figure}

Cellular structures have been a topic of interest for a range of bio-inspired materials because their geometries can be
quickly adapted to target desired mechanical properties.
In the work by Asbai-Ghoudan el at.\cite{Asbai2021} a structure for bone replacement is analyzed, and a key challenge identified when modeling such structures is mesh generation.
We construct models of similar geometry using the EB method to significantly simplify the meshing process.
The geometry of interest is defined by the gyroid function
\begin{equation}
    f(\xbold) = \cos(n \pi x)\sin(n \pi y)+\cos(n \pi y)\sin(n \pi z)+\cos(n \pi z)\sin(n \pi x) \,,
\end{equation}
evaluated in two dimensions along the $z=0.1$ plane where $f(\xbold) = 0$ with a period of $n=2$.
The surface is approximately thickened to width of $1/6$.
Flow is solved through a pipe of radius $r=0.95$ and length $8$, with a cut of the gyroid centered in the pipe at $x=4$.
A cut is made through the centerline with radius $r=0.175$.
The cut boundaries are smoothed (as in Devendran et al. \cite{Devendran2017}) to prevent singular solutions around sharp exterior corners.
Stokes flow through this structure is shown in \afigref{bone} for a grid refinement of $h=1/512$ in the flow direction.

\begin{figure}[!htbp]
  \centering
  \begin{subfigure}[t]{0.48\textwidth}
    \centering
    \includegraphics[width=0.9\linewidth]{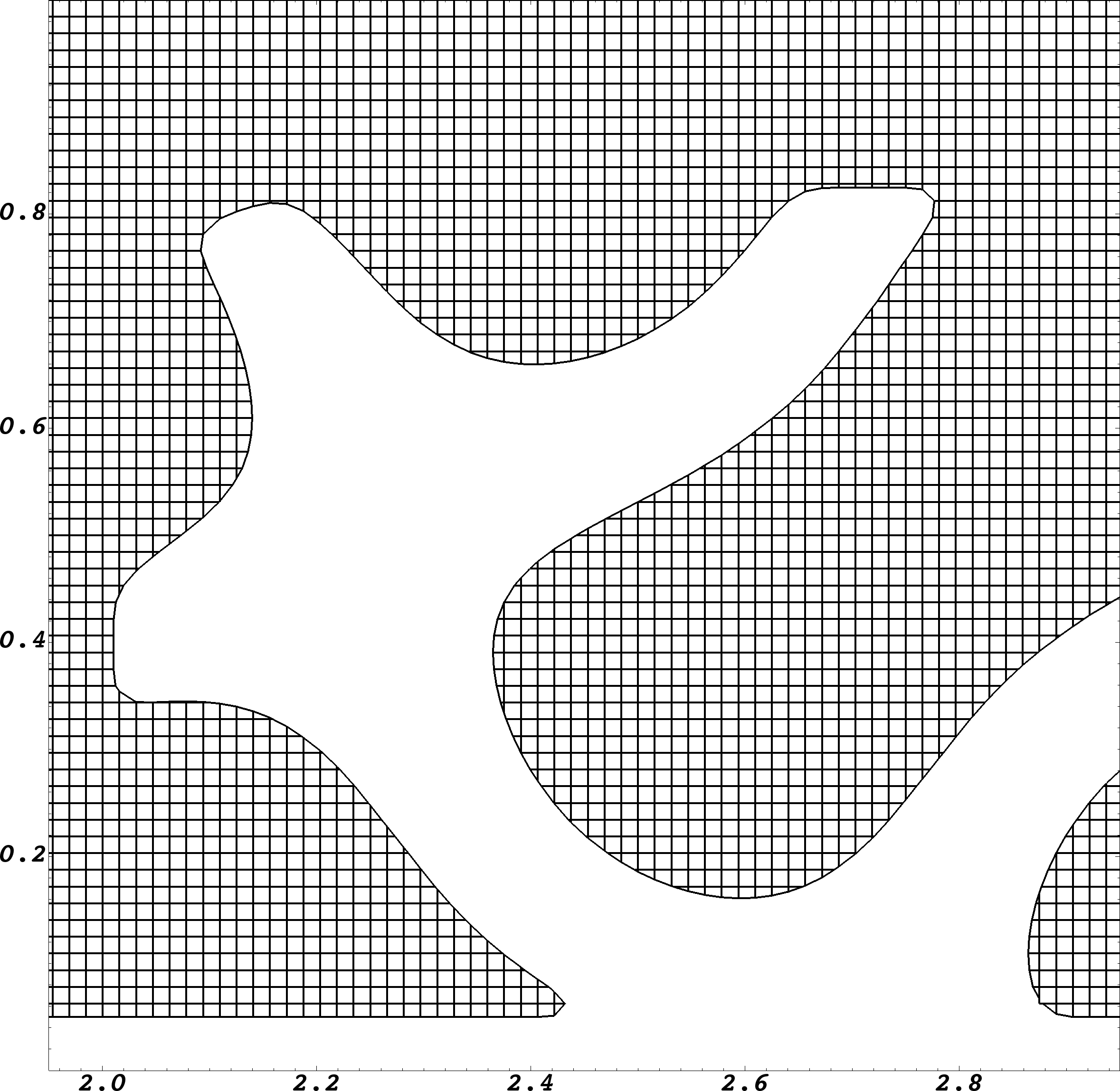}
    \caption{A close up of the EB grid.}
  \end{subfigure}
  \quad
  \begin{subfigure}[t]{0.48\textwidth}
    \centering
    \includegraphics[width=0.9\linewidth]{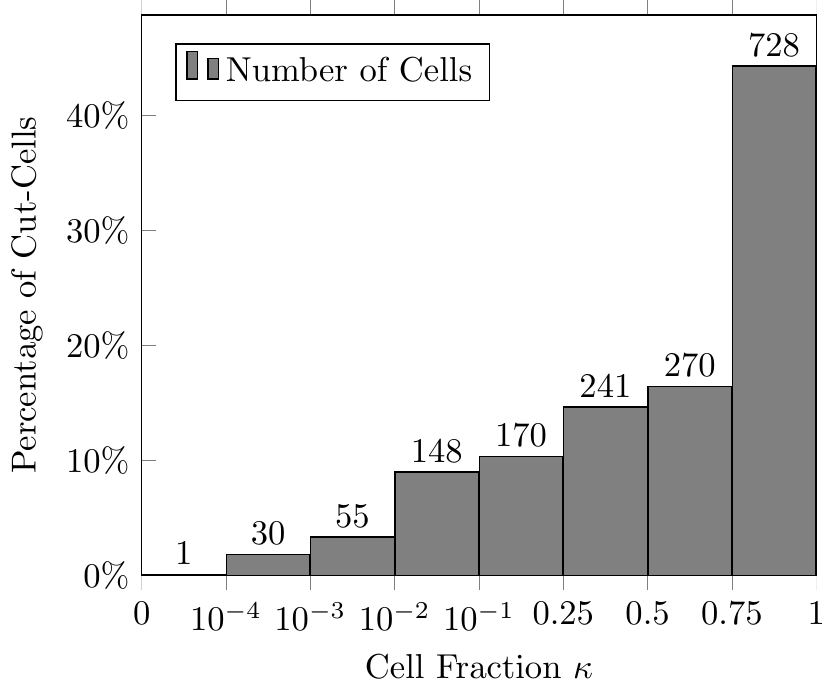}
    \caption{Distribution of the cut-cell sizes.}
  \end{subfigure}
  \caption{Bone scaffold grid.}
  \label{fig:boneGrid}
\end{figure}

The flow is found to be symmetric, and well-behaved along the cut-cell boundaries.
The EB grid for this case contains \num{26585} valid cells, of which \num{1643} are cut-cells with the smallest volume fraction of $\kappa$ = \num{5.933e-5}.
A close up view of this grid is shown in \afigref{boneGrid}(a), and the distribution of cut-cell sizes is shown in \afigref{boneGrid}(b).
In particular, the mesh generation requires trivial involvement, and because the geometry is specified analytically, adjustments such as gyroid size, position, and thickness are easily made.
This shows promise for more detailed analysis using the high-order EB method, where a large design space of geometries could be examined quickly without special considerations for mesh generation or solution stability.

%
\section{Conclusions and Future Work}\label{sec:conclusions}
In this work, we demonstrate that our EB algorithm capable of representing complex geometries is fourth-order accurate for the Stokes equations.
Additionally, our algorithm is shown to be stable when encountering small cells, without cell merging, redistribution, or grid remediation to avoid small cells.
These results demonstrate the feasibility of high-order EB methods, and provide confidence that the developed algorithm can be extended to solve engineering problems without making case specific considerations for mesh generation.

\subsection{Future Work}
The natural next step for this work is to tackle the incompressible Navier-Stokes equations.
This requires an approach for proper treatment of the non-linear advection term, and the ability to create stable upwind stencils for small cells.
We will also implement a geometric multigrid solver \cite{Martin1996, Zhang2012a} to improve the performance of our AMG-based solver.
Additionally, future work with this high-order EB method will include adaptive mesh refinement to improve solution accuracy in regions of interest \cite{Zhang2012a,Devendran2017}.


    
\end{document}